\documentclass[12pt,leqno]{article}
\usepackage{amsfonts}
\usepackage{amsmath, amsthm, amsfonts, amssymb, color}
\usepackage{mathrsfs}
\usepackage{url}
\usepackage{color}
\theoremstyle{definition}

\newcommand{\scr}[1]{\mathscr #1}
\definecolor{wco}{rgb}{0.5,0.2,0.3}

\numberwithin{equation}{section} \theoremstyle{remark}

\newcommand{\ua}{\uparrow}

\title{{\bf  Limit Theorems in  Wasserstein Distance  for Empirical Measures of Diffusion Processes on Riemannian Manifolds }\footnote{Supported in
 part by  NNSFC (11771326, 11831014, 11921001).} }
\author{
{\bf    Feng-Yu Wang$^{a),b)}$ \, \ Jie-Xiang Zhu$^{c),d)}$  }\\
\footnotesize{$^{a)}$ Center for Applied Mathematics, Tianjin University, Tianjin 300072, China }\\
 \footnotesize{ $^{b)}$ Department of Mathematics,
Swansea University,
Bay Campus,
Swansea,
SA1 8EN, United Kingdom}\\
\footnotesize{$^{c)}$ School of Mathematical Sciences, Fudan University, Shanghai 200433, China }\\
\footnotesize{$^{d)}$ Institut de Math\'ematiques de Toulouse, Universit\'e de Toulouse -- Paul-Sabatier, F-31062 Toulouse, France }\\
\footnotesize{  wangfy@tju.edu.cn, F.-Y.Wang@swansea.ac.u; 15110840006@fudan.edu.cn, jiexiangzhu7@gmail.com } }

\begin{document}
\allowdisplaybreaks
\def\R{\mathbb R}  \def\ff{\frac} \def\ss{\sqrt} \def\B{\mathbf
B}\def\TO{\mathbb T}
\def\I{\mathbb I_{\pp M}}\def\p<{\preceq}
\def\N{\mathbb N} \def\kk{\kappa} \def\m{{\bf m}}
\def\ee{\varepsilon}\def\ddd{D^*}
\def\dd{\delta} \def\DD{\Delta} \def\vv{\varepsilon} \def\rr{\rho}
\def\<{\langle} \def\>{\rangle} \def\GG{\Gamma} \def\gg{\gamma}
  \def\nn{\nabla} \def\pp{\partial} \def\E{\mathbb E}
\def\d{\text{\rm{d}}} \def\bb{\beta} \def\aa{\alpha} \def\D{\scr D}
  \def\si{\sigma} \def\ess{\text{\rm{ess}}}
\def\beg{\begin} \def\beq{\begin{equation}}  \def\F{\scr F}
\def\Ric{{\rm Ric}} \def\Hess{\text{\rm{Hess}}}
\def\e{\text{\rm{e}}} \def\ua{\underline a} \def\OO{\Omega}  \def\oo{\omega}
 \def\tt{\tilde}
\def\cut{\text{\rm{cut}}} \def\P{\mathbb P} \def\ifn{I_n(f^{\bigotimes n})}
\def\C{\scr C}      \def\aaa{\mathbf{r}}     \def\r{r}
\def\gap{\text{\rm{gap}}} \def\prr{\pi_{{\bf m},\varrho}}  \def\r{\mathbf r}
\def\Z{\mathbb Z} \def\vrr{\varrho} \def\ll{\lambda}
\def\L{\scr L}\def\Tt{\tt} \def\TT{\tt}\def\II{\mathbb I}
\def\i{{\rm in}}\def\Sect{{\rm Sect}}  \def\H{\mathbb H}
\def\M{\scr M}\def\Q{\mathbb Q} \def\texto{\text{o}} \def\LL{\Lambda}
\def\Rank{{\rm Rank}} \def\B{\scr B} \def\i{{\rm i}} \def\HR{\hat{\R}^d}
\def\to{\rightarrow}\def\l{\ell}\def\iint{\int}
\def\EE{\scr E}\def\Cut{{\rm Cut}}\def\W{\mathbb W}
\def\A{\scr A} \def\Lip{{\rm Lip}}\def\S{\mathbb S}
\def\BB{\scr B}\def\Ent{{\rm Ent}} \def\i{{\rm i}}\def\itparallel{{\it\parallel}}
\def\g{{\mathbf g}}\def\Sect{{\mathcal Sec}}\def\T{\mathcal T}\def\V{{\bf V}}
\def\PP{{\bf P}}\def\HL{{\bf L}}\def\Id{{\rm Id}}\def\f{{\bf f}}\def\cut{{\rm cut}}

\def\BL{\scr A}

\maketitle

\begin{abstract} Let $(M,\rr)$ be a   connected compact Riemannian manifold without boundary or with    a  convex boundary $\pp M$, let $V\in C^2(M)$ such that $\mu(\d x):=\e^{V(x)}\d x$ is a probability measure, where $\d x$ is the volume measure. Let
   $\{\ll_i\}_{i\ge 1} $ be all non-trivial eigenvalues of $-L$  with Neumann boundary condition if $\pp M\ne \emptyset$, where $L := \DD + \nn V$ for $\DD$ being the Laplace-Beltrami operator on $M$. Then  the  empirical measures $\{\mu_{t}\}_{t>0}$ of the diffusion process generated by $L$ (with reflecting boundary if $\pp M\ne\emptyset$) satisfy
$$ \lim_{t\to \infty}  \big\{t \E^x [\W_2(\mu_{t},\mu)^2]\big\}= \sum_{i=1}^\infty\ff 2 {\ll_i^2}\ \text{ uniformly\ in\ }  x\in M,$$ where  $\E^x$ is the expectation   for the diffusion process starting at point $x$, and
    $\W_2$ is the  $L^2$-Wasserstein  distance induced by the Riemannian metric. The limit is finite if and only if $d\le 3$, and in this case    we   derive     $$\lim_{t\to\infty} \sup_{x\in M} \Big|\P^x(t  \W_2(\mu_{t},\mu)^2<a)- \P\Big( \sum_{k=1}^\infty \ff{2\xi_k^2}{\ll_k^2}<a\Big)\Big|=0, \ \ a\ge 0,$$   where $\P^x$ is the probability with respect to $\E^x$, and $\{\xi_k\}_{k\ge 1}$ are i.i.d. standard Gaussian random variables.
    Moreover,
   $\E^x[\W_2(\mu_{t},\mu)^2]\sim t^{-\ff 2 {d-2}}$ for $d\ge 5$,  and when  $d=4$ we have  $\E^x[\W_2(\mu_{t},\mu)^2] \le c    t^{-1}\log t$  for some constant $c>0$ and large $t$ while the same type lower bound estimate holds for $M=\TO^4$. Finally,  we establish  the long-time large deviation principle for $\{\W_2(\mu_t,\mu)^2\}_{t\ge 0}$   with  a good rate function given by  the information with respect to $\mu$.
 \end{abstract} \noindent
 \vskip 2cm

\section{Introduction and Main results}

The diffusion processes (for instance, the Brownian motion) on Riemannian manifolds have intrinsic link  to  properties (for instances, curvature, dimension, spectrum) of the infinitesimal generator, see, for instances, the monographs \cite{BIL,W14}  and references within. In this paper,  we  characterize   the long time behaviour of empirical measures for   diffusion processes by using eigenvalues of the generator.

Let $M$ be a $d$-dimensional    connected  complete Riemannian manifold   possibly with a smooth  boundary $\pp M$.  Let $V\in C^2(M)$ such that $\mu(\d x)=\e^{V}\mu_M(dx)$ is a probability measure on $M$, where $\mu_M$ is the Riemannian volume measure on $M$. Then  the   (reflecting, if $\pp M\ne \emptyset$) diffusion process $X_t$ generated by $L:=\DD+\nn V$ on $M$ is reversible; i.e. the associated diffusion semigroup $\{P_t\}_{t\ge 0}$ is symmetric in $L^2(\mu)$, where
 $$P_tf(x):= \E^x f(X_t),\ \ t\ge 0, f\in \B_b(M). $$
  Here, $\E^x$ is  the expectation  taken for the diffusion process $\{X_t\}_{t\ge 0} $ with
 $X_0=x,$ and we will use $\P^x$ to denote the associated probability measure.
 In general, for any probability measure $\nu$ on $M$, let $\E^\nu$ and $\P^\nu$ be the expectation and probability taken for the diffusion process with initial distribution $\nu$.

When the diffusion process generated by $L$ is exponentially ergodic,  it is in particular the case  when $M$ is compact,  the empirical measure  $$\mu_t:=\ff 1 t \int_0^t \dd_{X_s} \d s,\ \ t>0$$ converges weakly to $\mu$ as $t\to\infty$. More precisely,
  for any non-constant $f\in C_b(M)$,  we have the law of large number
$$\lim_{t\to\infty} \mu_t(f)= \mu(f)\ \text{a.s.}$$
as well as the central limit theorem
$$\ss{t} \big\{\mu_t(f) -\mu(f)\big\}\to N(0,\dd(f))\ \text{in\ law\ as\ }t\to\infty,$$ where $\dd(f):=\lim_{t\to\infty}   t \E|\mu_t(f)-\mu(f)|^2\in (0,\infty)$ exists, and $N(0,\dd(f))$ is the centered normal distribution with variance $\dd(f)$.  Consequently,
  the average  additive functional $\mu_t(f)$ converges to $\mu(f)$ in $L^2(\P)$ with rate  $t^{-\ff 1 2}$, which is universal and has nothing to do with specific properties of $M$ and $L$.
 See, for instance  \cite{Kulik},  for historical remarks and  more results concerning limit theorems on additive functionals of Markov processes.

 On the other hand, since the Wasserstein distance $\W_2$ induced by the Riemannian distance $\rr$ on $M$  is  associated with  a natural Riemannian structure on the space of probability measures, see e.g. \cite{Otto}, the asymptotic behaviors of $\W_2(\mu_t,\mu)$
  should reflect intrinsic  properties of $M$ and  $L$.  Indeed, as shown in Theorem \ref{T1.1} below, the  long time behavior of  $\W_2(\mu_t,\mu)^2$ depends on the dimension of $M$ and all eigenvalues of $L$, this is essentially different from that of the additive functional
  $\mu_t(f)$ introduced above.

Let $\scr P$ be the set of all probability measures on $M$, and let   $\rr$ be the Riemiannian distance on $M$.  For any $p\ge 1$,  the $L^p$-Wasserstein distance $\W_p$ is defined by
 $$\W_p(\mu_1,\mu_2):= \inf_{\pi\in \C(\mu_1,\mu_2)} \bigg(\int_{M\times M} \rr(x,y)^p \pi(\d x,\d y) \bigg)^{\ff 1 p},\ \ \mu_1,\mu_2\in \scr P,$$
 where $\C(\mu_1,\mu_2)$ is the set of all probability measures on $M\times M$ with marginal distributions $\mu_1$ and $\mu_2$. A measure $\pi\in \C(\mu_1,\mu_2)$ is called a coupling of $\mu_1$ and $\mu_2$.

 In this paper, we aim to characterize the long time behavior  of  $ \W_2(\mu_t,\mu)^2$.
When $M$ is compact,  we will prove the large deviation principle with rate function
$$I(r):=\inf\{I_\mu(\nu):\ \nu\in\scr P, \W_2(\nu,\mu)\ge r\},\ \ r\ge 0,$$ where $I_\mu$ is the information with respect to $\mu$; i.e.
$$I_\mu(\nu):=\beg{cases} \mu(|\nn f^{\ff 1 2} |^2), &\text{if}\ \nu=f\mu, f^{\ff 1 2}\in W^{2,1}(\mu);\\
\infty, &\text{otherwise}.\end{cases} $$ Here, $W^{2,1}(\mu)$ is the closure of $C^\infty(M)$ under the Sobolev norm
$$\|h\|_{2,1}:=\ss{\mu(h^2+|\nn h|^2)}.$$
By convention, we set $\inf\emptyset =\infty$, so that $I(r)=\infty$ for $r>r_0$, where since $\rr$ is bounded,
$$r_0:= \sup_{\nu}\W_2(\nu,\mu)^2 = \sup_{x\in M} \mu(\rr(x,\cdot)^2)<\infty.$$

It is well known that when $M$ is compact, $L$ has purely discrete spectrum, and all non-trivial eigenvalues $\{\ll_i\}_{i\ge 1}$ of $-L$ listed in the increasing order counting multiplicities satisfy (see for instance \cite{Chavel})
 \beq\label{*2} \kk^{-1} i^{\ff 2 d} \le \ll_i\le \kk i^{\ff 2 d},\ \ i\ge 1\end{equation}
for some constant $\kk>1$. Our first result is the following.

\beg{thm}\label{T1.1}  Let $M$ be compact.
\beg{enumerate} \item[$(1)$] If $\pp M$ is empty or convex, then the following limit  formula holds uniformly in $x\in M$:
\beq\label{ULT}   \lim_{t\to \infty}     \Big\{t \E^x [\W_2(\mu_{t},\mu)^2]\Big\}=    \sum_{i=1}^\infty\ff 2 {\ll_i^2}.\end{equation}
In general, there exists a constant $c\in (0,1]$ such that
\beq\label{UEST} \beg{split}  &c\sum_{i=1}^\infty\ff 2 {\ll_i^2} \le \liminf_{t\to \infty}  \inf_{x\in M}    \Big\{t \E^x [\W_2(\mu_{t},\mu)^2]\Big\}\\
&\le \limsup_{t\to \infty}  \sup_{x\in M}    \Big\{t \E^x [\W_2(\mu_{t},\mu)^2]\Big\}
\le   \sum_{i=1}^\infty\ff 2 {\ll_i^2},\end{split} \end{equation}
\item[$(2)$]  $\{\W_2(\mu_t,\mu)^2 \}_{t\ge 0} $ satisfies the   uniform large deviation principle with good rate function $I$; that is, $\{I\le \aa\}$ is a compact subset of $[0,\infty)$ for any $\aa\in [0,\infty)$, and
\beg{align*} &-\inf_{r\in A^\circ} I(r)\le \liminf_{t\to\infty} \ff 1 t \log  \inf_{x\in M} \P^x(\W_2(\mu_t,\mu)^2\in A^\circ) \\
&\le \limsup_{t\to\infty} \ff 1 t \log  \sup_{x\in M} \P^x(\W_2(\mu_t,\mu)^2\in \bar A)  \le - \inf_{r\in \bar A} I(r),\ \ A\subset [0,\infty),\end{align*} where $A^\circ$ and $\bar A$ denote the interior and the closure of $A$ respectively.
\item[$(3)$] If $d\le 3$ and $\pp M$ is either empty or convex, then
 \beq\label{CLT'} \lim_{t\to \infty}\sup_{x\in M} \big|\P^x\big(t  \W_2(\mu_{t},\mu)^2 <a\big)-\nu_0\big((-\infty,a)\big)\big|=0,\ \ a\in \R,\end{equation} where $\nu_0$ is the distribution
  of\, $\Xi_0:= \sum_{k=1}^\infty \ff{2\xi_k^2}{\ll_k}$ for a sequence of i.i.d. standard Gaussian random variables $\{\xi_k\}_{k\ge 1}.$ \end{enumerate}
  \end{thm}

In Theorem \ref{T1.1}(3) we only consider $d\le 3$, since  Theorem \ref{T1.1}(1) and \eqref{*2} yield
\beq\label{*DP} \liminf_{t\to\infty} \Big\{t \inf_{x\in M} \E^x [\W_2(\mu_{t},\mu)^2] \Big\}=\infty,\ \ d\ge 4.\end{equation}
So, for $d\ge 4$  the convergence  of  $\E^x [\W_2(\mu_{t},\mu)^2]$ is slower than $t^{-1}$.
In the next result we     present    two-sided estimates on the convergence rate of $\E[\W_2(\mu_t,\mu)^2]$ for $d\ge 4$.

\beg{thm}\label{T1.1'} Let $M$ be compact with $d\ge 4$. \beg{enumerate} \item[$(1)$] There exists a constant $c>0$ such that for any $t\ge 1$,
$$  \sup_{x\in M}  \E^x[\W_2(\mu_t,\mu)^2] \le \beg{cases} c t^{-1} \log (1+t), &\text{if}\ d=4,\\
c t^{-\ff 2{d-2}}, &\text{if}\ d\ge 5.\end{cases}$$
 \item[$(2)$] On the other hand, there exists a constant $c'>0$ such that
$$ \inf_{x\in M}  \E^x[\W_2(\mu_t,\mu)^2] \ge  \inf_{x\in M}  \{\E^x[\W_1(\mu_t,\mu)] \}^2 \ge  c' t^{-\ff 2 {d-2}},\ \ t \ge 1.$$
\end{enumerate}
\end{thm}
 Theorem \ref{T1.1'} implies that when $d\ge 5$ we have  $\E [\W_2(\mu_t,\mu)^2]\sim t^{-\ff 2 {d-2}}$ for large $t$.
 Due to \eqref{*DP},    we hope that $\E [\W_2(\mu_t,\mu)^2] \sim t^{-1}\log t$ holds for $d=4$, i.e. the order in the upper bound  estimate is  sharp.
Although in general this is not yet proved in the paper, it is true    for    $M=\TO^4 (=\R^4 \backslash (2\pi \Z^4))$   and $V=0$  according to  the following result and  $\{\E\W_1(\mu_t,\mu)\}^2\le \E [\W_2(\mu_t,\mu)^2]. $

\beg{thm}\label{TN} Let $M=\TO^4$ and $V=0$. Then there exists a constant $c>0$ such that
$$\inf_{x\in M}\{ \E^x \W_1(\mu_t,\mu)\}^2 > c t^{-1} \log t,\ \ t \ge 1$$ \end{thm}

 \

To conclude this section,
 we  compare the convergence rate of  $\W_2(\mu_t,\mu)$  with that of $\W_2(\bar \mu_n,\mu)^2$ investigated  in \cite{[2], AMB,  [9], BLG14, [11], [21], [25]},
 where  $$\bar\mu_n:=\ff 1 n\sum_{i=1}^n \dd_{X_i},\ \ n\in\mathbb N$$   is  the empirical measure of  i.i.d. random variables   $\{X_i\}_{i\ge 1}$ with common distribution $\mu$.
 In particular,  for $\mu$ being the uniform distribution on a bounded domain in $\R^d$, we have
 $$
 \E[ \W_2(\bar \mu_n,\mu)^2] \sim  \beg{cases} n^{-2/d}\ &\text{if\ }d\ge 3,\\
  n^{-1}\log n, \ &\text{if\ }d=2,\\
 n^{-1}\ &\text{if}\ d=1,
\end{cases}
  $$
 where    the   assertion for $d=2$ is known as Ajtai-Koml\'os-Tusn\'ady (AKT) optimal matching theorem \cite{[2]}, and $a_n\sim b_n$ means that $c_1a_n\le b_n\le c_2 a_n$ holds for some constants $c_2\ge c_1>0$ and large $n$. Moreover, the empirical measure   for a discrete time Markov chain on a bounded domain has been  investigated in \cite{BLG14}.

Combining this Theorems \ref{T1.1}-\ref{TN} for $t=n$, we see that   in the sense of    empirical measures,  diffusion processes  converge  faster  than   i.i.d. samples:  
  the empirical measures of $d$-dimensional diffusion processes   behave as  those of   $(d-2)\lor 1$-dimensional i.i.d. samples. 
In particular, in the present setting the feature of AKT optimal matching theorem appears to dimension $d=4$ rather than $d=2$. However, unlike in the i.i.d. case for which  Ambrosio-Stra-Trevisan   \cite{AMB} derived
 the exact value of $\lim_{n\to\infty} \ff n {\log n} \E[\W_2(\bar\mu_n,\mu)^2] $ for the uniform distribution $\mu$ on $\TO^2$,  in the present setting
  the exact value of $ \lim_{t\to \infty}\ff t{\log t}\E [\W_2(\mu_{t},\mu)^2]$ is   unknown for the uniform distribution $\mu$ on $\TO^4$.   This could be a challenging problem.

Since $\mu_t$ is singular with respect to $\mu$, it is hard to estimate $\W_2(\mu_t,\mu)$ using analytic methods.  So, as in \cite{AMB0},  we first investigate   the modified empirical measures
$$\mu_{t,r}  := \mu_t P_r= \ff 1 t\int_0^t \{\dd_{X_s}P_r\} \d s, \ \ \ t>0,\ r>0,$$
where  for a  probability measure $\nu$ on $M$,   $\nu P_r$ denotes the distribution of $X_r$ with $X_0$ having law $\nu$.
Note that $\lim_{r\to 0} \W_2(\mu_{t,r}, \mu_t)= 0$, see \eqref{LK1} below for an estimate of the convergence rate. The main results of the paper have been extended in \cite{WW} 
to subordinated diffusion processes on compact manifolds, see also \cite{W1,W2,W3,W4} for further development on the empirical measures for killed diffusion processes, SDEs and SPDEs.

The remainder of the paper is organized as follows. In Section 2,  we investigate the long time behavior of the modified empirical measures $\mu_{t,r}$ for $r>0$, where $M$ might be non-compact.
We then prove  Theorems \ref{T1.1},   \ref{T1.1'}   and \ref{TN} in Sections 3,  4  and 5 respectively.

 \section{Asymptotics for modified empirical measures}

In this part, we allow $M$ to be non-compact, but assume that $L$ satisfies the curvature condition
  \beq\label{CV} \Ric_V:= \Ric-\Hess_V\ge -K\end{equation} for some constant $K\ge 0$, where $\Ric$ is the Ricci curvature on $M$ and $\Hess_V$ is the Hessian tensor of $V$.   This condition means that
  $\Ric_V(X,X)\ge -K|X|^2$ for all $X\in TM$, the tangent bundle of $M$.

When $\pp M\ne \emptyset$, let $N$ be the inward unit normal vector field of $\pp M$. We call $\pp M$ convex, if its  second fundamental form
$\I$ is nonnegative; i.e.
$$\I(X,X):=-\<X,\nn_XN\>\ge 0,\ \ X\in T\pp M,$$ where $T\pp M$ is the tangent bundle of the boundary $\pp M$.
In general, for a function $g$ on $\pp M$, we write
$\I\ge g$ if
\beq\label{GH} \I(X,X)\ge g(x)|X|^2,\ \ x\in \pp M, X\in T_x \pp M.\end{equation}  We call $\pp M$ convex on a set $D\subset M$, if \eqref{GH} holds for some function $g$ which is non-negative on $D\cap \pp M$.

For any $q\ge p\ge 1$, let $\|\cdot\|_{p\to q}$ be the operator norm from $L^p(\mu)$ to $L^q(\mu)$. We will need the following assumptions.

\beg{enumerate}   \item[{\bf (A1)}]    $P_t$ is ultracontractive, i.e.  $\|P_t\|_{1\to \infty}:=\sup\limits_{\mu(|f|)\le 1} \|P_t f\|_\infty <\infty, \ t>0.$
\item[{\bf (A2)}]  $\eqref{CV}$ holds for some constant $K\ge 0$, and there exists a compact set $D\subset M$ such that either $D^c\cap \pp M=\emptyset$ or $\pp M$ is convex on $D^c$.
\end{enumerate}

 Obviously, {\bf (A1)} and  {\bf (A2)} hold   if $M$ is compact.
When $M$ is non-compact satisfying condition {\bf (A2)}, by  \cite[Theorem 3.5.5]{W14},    {\bf (A1)} holds if and only  if  $\|P_t\e^{\ll \rr_o(\cdot)^2}\|_\infty<\infty $  for any $\ll,t>0$,
where  $\rr_o:=\rr(o,\cdot)$ is the distance function to a fixed point $o\in M$, see  \cite[Corollary 2.5]{RW03} for concrete examples with $\|P_t\e^{\ll \rr_o(\cdot)^2}\|_\infty<\infty$.
 See also \cite[Proposition 4.1]{W09} for  examples satisfying assumption  {\bf (A1)} when $\Ric_V$ is unbounded from below.

{\bf (A1)} implies that the spectrum of $L$ (with Neumann boundary condition if $\pp M\ne \emptyset$)  is purely discrete.
Since $M$ is connected, in this case $L$ has a spectral gap, i.e. $0$ is a simple  isolated eigenvalue of $L$.
Let $\{\ll_i\}_{i\ge 1} $ be all non-trivial eigenvalues of $-L$  listed in the increasing order including multiplicities. By the concentration of $\mu$ implied by the ultracontractivity condition {\bf (A1)}, we have
\beq\label{MM} \int_{M\times M} \e^{\ll \rr^2}\d(\mu\times\mu)<\infty,\ \ \ll>0.\end{equation}
Indeed, according to \cite{DS, Gross} (see for instance \cite[Theorem 1.1]{RW03}), {\bf (A1)} implies that for some $\bb: (0,\infty)\to (0,\infty)$,
$$\mu(f^2\log f^2)\le r \mu(|\nn f|^2)+\bb(r),\ \ r>0, f\in C_b^1(M), \mu(f^2)=1,$$
which then ensures \eqref{MM} by \cite[Corollary 6.3]{RW03} or \cite{AMS}.

 For any $r>0$, let $\nu_r$ be the distribution of $$\Xi_r:=\sum_{k=1}^\infty \ff{2\xi_k^2} {\ll_k^2\e^{2\ll_k r}},$$ where  $\{\xi_k\}_{k\ge 1}$ are i.i.d. standard Gaussian random variables.

\beg{thm}\label{T1.2}   Assume    {\bf (A1)} and let $r>0$. Then
 \beq\label{E1'}  \limsup_{t\to \infty} \sup_{x\in M} \Big\{t \E^x [\W_2(\mu_{t,r},\mu)^2]\Big\}\le   \sum_{i=1}^\infty\ff 2 {\ll_i^2\e^{2r\ll_i}}<\infty.\end{equation}
 If moreover {\bf (A2)} holds, then
\beq\label{E2'} \lim_{t\to \infty} \sup_{x\in M} \bigg|t \E^x [\W_2(\mu_{t,r},\mu)^2]- \sum_{i=1}^\infty\ff 2 {\ll_i^2\e^{2r\ll_i}} \bigg| =0,\end{equation}
 \beq\label{CLT} \lim_{t\to \infty}\sup_{x\in M} \big|\P^x\big(t  \W_2(\mu_{t,r},\mu)^2 <a\big)-\nu_r\big((-\infty,a)\big)\big|=0,\ \ a\in \R.\end{equation}
\end{thm}
\paragraph{Remark 2.1.} Consider the measure
$$\mu_{sp}:= \sum_{i=1}^\infty \ff 1 {\ll_i^2} \dd_{\ll_i},$$ whose support consists of  all non-trivial eigenvalues of $L$. Then \eqref{E2'}  implies
$$\int_0^\infty \e^{-2r s}\mu_{sp}(\d s) =\lim_{t\to\infty} \Big\{t \E^\nu [\W_2(\mu_{t,r},\mu)^2]\Big\},\ \ r> 0$$ for any probability measure $\nu$ on $M$.
This gives a  probabilistic representation for the Laplace transform of $\mu_{sp}$, and hence determines all eigenvalues and multiplicities for $L$.

\

To investigate the long time behavior of $ \E[\W_2(\mu_{t},\mu)^2]$, i.e. $ \E[\W_2(\mu_{t,r},\mu)^2]$ with $r=0$, one may consider the limit of formula \eqref{E2'} when  $r\downarrow 0$.
\beg{cor}\label{C1.2} If $M$ is compact, then:
\beg{enumerate}\item[$(1)$] For $d\le 3$,
$$ \lim_{r\downarrow 0}\lim_{t\to \infty}   \Big\{t \E^x [\W_2(\mu_{t,r},\mu)^2]\Big\}=   \sum_{i=1}^\infty\ff 2 {\ll_i^2}<\infty\ \text{uniformly\ in\ }x\in M.$$
\item[$(2)$] For $d=4,$
$$\lim_{r\downarrow 0} \lim_{t\to \infty}   \ff{\log \log  \{t \E^x [\W_2(\mu_{t,r},\mu)^2]\}}{\log\log r^{-1}}=1\ \text{uniformly\ in\ }x\in M.$$
\item[$(3)$] For $d\ge 5,$
$$\lim_{r\downarrow 0} \lim_{t\to \infty}   \ff{\log    \{t \E^x [\W_2(\mu_{t,r},\mu)^2]\}}{\log  r^{-1}}=\ff{d-4}2\ \text{uniformly\ in\ }x\in M.$$\end{enumerate}\end{cor}

  In the following two  subsections, we   investigate the upper and lower bound estimates on
  $\E[\W_2(\mu_{t,r},\mu)^2]$   respectively,  which then  lead to proofs of Theorem \ref{T1.2} and Corollary \ref{C1.2} in the last subsection.

\subsection{Upper bound estimate}

We first   estimate $\W_2(\mu_1,\mu_2)$  in terms of the energy for the difference of the density functions of $\mu_1$ and $\mu_2$ with respect to $\mu$.
Let   $\D(L)$ be  the domain of the generator $L$ in $L^2(\mu)$, with Neumann boundary condition if $\pp M\ne \emptyset$.
Then for any function
$$g\in L^2_0(\mu):=\{g \in L^2(\mu), \mu(g) = 0 \},$$ we have
\beq\label{LLL}  (-L)^{-1}g =\int_0^\infty P_sg\,\d s\in \D(L),\ L(L^{-1})g =g. \end{equation}  Since $M$ is complete and $\mu$ is finite,
we have $\D(L)\subset \D((-L)^{\ff 1 2})= H^{1,2}(\mu)=W^{1,2}(\mu)$, where $H^{1,2}(\mu)$ is the completion of $C_0^\infty(M)$ under the Sobolev norm
$$\|f\|_{1,2}:= \ss{\mu(f^2) +\mu(|\nn f|^2)},$$ and $W^{1,2}(\mu)$ is the class of all weakly differentiable functions $f$ on $M$ such that $|f|+|\nn f|\in L^2(\mu).$ In particular,   $L^{-1}g\in W^{1,2}(\mu)$ for $g\in L_0^2(\mu)$. The following lemma  is  essentially due to \cite[Proposition 2.3]{AMB} where the case with compact $M$ and  $V=0$ is concerned, but its proof works also for the present setting.

\beg{lem} \label{L2.1} Let $f_0,f_1\in L^2(\mu)$ be probability density functions with respect to $\mu$.   Then
 $$\W_2(f_0\mu, f_1\mu)^2\le \int_M \ff{|\nn L^{-1}(f_1-f_0)|^2}{\M(f_0,f_1)}\d\mu,$$
 where   $\M(a,b):= \ff{a-b}{\log a-\log b}$ for $a,b>0$, and  $\M(a,b):=0$ if one of $a$ and $b$ is zero.
 \end{lem}

 \beg{proof} Let ${\rm Lip}_b(M)$ be  the set of bounded Lipschitz continuous functions  on $M$.    Consider the Hamilton-Jacobi semigroup $(Q_t)_{t>0}$ on ${\rm Lip}_b(M)$:
 $$ Q_t \phi:= \inf_{x\in M} \Big\{\phi(x)+ \ff 1 {2 t} \rr(x,\cdot)^2\Big\},\ \ t>0, \phi\in {\rm Lip}_b(M).$$
 Then for any $\phi\in {\rm Lip}_b (M)$, $Q_0\phi:= \lim_{t\downarrow 0} Q_t\phi=\phi$, $\|\nn Q_t\phi\|_\infty$ is locally bounded in $t\ge 0$, and $Q_t \phi$ solves the  Hamilton-Jacobi equation
 \beq\label{HK0} \ff{\d}{\d t} Q_t \phi= -\ff 1 2 |\nn Q_t\phi|^2,\ \ t>0.\end{equation}
 In a more general setting of metric spaces, one has  $\ff{\d}{\d t} Q_t \phi\le -\ff 1 2 |\nn Q_t\phi|^2\ \mu$-a.e.,
  where the equality holds for length spaces which include the present framework, see e.g. \cite{AMB0, AMB}.

Letting $\mu_i= f_i \mu, i=0,1,$    the Kantorovich dual formula implies
 \beq\label{KT} \ff 1 2 \W_2(\mu_0,\mu_1)^2= \sup_{\phi\in {\rm Lip}_b(M)} \big\{\mu_1(Q_1\phi)- \mu_0(\phi)\big\}.\end{equation}
 Let $f_s= (1-s)f_0+s f_1, s\in [0,1].$ By \eqref{MM} and the boundedness of $\|\nn Q_t\phi\|_\infty$ in $t\in [0,1]$, we deduce from \eqref{HK0} that
 \beq\label{KT'}\ff{\d}{\d s}  \int_M f_sQ_s\phi\d\mu  =\int_M \Big\{-\ff 1 2 |\nn Q_s\phi|^2 f_s+(Q_s\phi) (f_1-f_0)\Big\}\d\mu,\ \ s\in (0,1].\end{equation}
 Moreover,   \eqref{LLL} implies  $f:= L^{-1} (f_0-f_1)\in \D(L)$. Then by \eqref{KT'} and using
 the integration by parts formula,  for any $\phi\in {\rm Lip}_b(M)$ we have
 \beg{align*} &\mu_1(Q_1\phi)- \mu_0(\phi)=\int_M \big\{f_1Q_1\phi-f_0\phi\big\}\d\mu =\int_0^1 \bigg(\ff{\d}{\d s} \int_M f_s Q_s \phi\, \d\mu \bigg)\d s\\
 &=  \int_0^1 \d s \int_M \Big\{-\ff 1 2 |\nn Q_s \phi|^2 f_s + (Q_s \phi)  (f_1-f_0) \Big\}\d\mu  \\
 & =\int_0^1 \d s \int_M \Big\{-\ff 1 2 |\nn Q_s \phi|^2 f_s - (Q_s \phi)  L f  \Big\}\d\mu \\
 &= \int_0^1 \d s \int_M \Big\{-\ff 1 2 |\nn Q_s \phi|^2 f_s + \<\nn f, \nn  Q_s \phi \> \Big\}\d\mu
  \le \ff 1 2  \int_0^1\d s  \int_M \ff{|\nn f|^2}{f_s}\d \mu \\
  &= \ff 1 2\int_M|\nn f|^2 \d\mu \int_0^1  \ff {\d s} {(1-s)f_0+sf_1} = \ff 1 2 \int_M \ff{|\nn f|^2}{\M(f_0,f_1)}\d\mu.\end{align*}
 Combining this with \eqref{KT}, we finish the proof.
 \end{proof}

 To estimate $ \W_2(\mu_{t,r},\mu)^2$ using   Lemma \ref{L2.1},    we need to figure out  the density function $f_{t,r} $ of $\mu_{t,r}$ with respect to $\mu$,
 i.e. $f_{t,r}$ is a nonnegative function such that $\mu_{t,r}(A)= \int_A f_{t,r}\d\mu$  for any measurable set  $A\subset M$.
Obviously, letting $p_t(x,y)$   be the heat kernel of $P_t$ with respect to $\mu$, i.e.
$$P_tf(x)=\int_M p_t(x,y)f(y)\mu(\d y),\ \ t>0, x\in M, f\in \B_b(M),$$
 we have
 \beq\label{FT} f_{t,r}:= \ff 1 t\int_0^t p_r(X_s,\cdot)\d s,\ \ t>0.\end{equation}
On the other hand, the assumption   {\bf (A1)} implies
\beq\label{HTT} \sup_{x,y\in M} p_t(x,y)=\|P_t\|_{1\to\infty}<\infty,\ \ t>0,\end{equation}
\beq\label{HK} p_t(x,y)=1+ \sum_{i=1}^\infty \e^{-\ll_i t} \phi_i(x)\phi_i(y),\ \ t>0, x,y\in M,\end{equation}
where  $\{\phi_i\}_{i\ge 1}$ are unit (Neumann if $\pp M\ne \emptyset$) eigenfunctions of  $-L$ with eigenvalues $\{\ll_i\}_{i\ge 1}.$
In particular, \eqref{HK} implies
  \beq\label{NHK} f_{t,r}(y) -1 = \ff 1 {\ss t} \sum_{i=1}^{\infty} \e^{-\ll_i r} \psi_i(t) \phi_i(y),\ \ \psi_i(t):= \ff 1 {\ss t} \int_0^t \phi_i(X_s)\d s, \ i \in\mathbb N.\end{equation}
The assumption   {\bf (A1)}  also implies the following  Poincar\'e inequality,
\beq\label{SPP} \|P_tf\|_2\le \e^{-\ll_1 t} \|f\|_2,\ \ t \ge 0, f\in L^2_0(\mu).\end{equation}
 Since  $P_t$ is contractive in $L^p(\mu)$ for any $p\in [1,\infty]$,    \eqref{HTT} and \eqref{SPP}    yield
 \beq\label{ULT-2} \|P_tf\|_p\le c \e^{-\ll_1 t}\|f\|_p,\ \ t\ge 0, p\in [1,\infty], f\in L^p_0(\mu) \end{equation}  for some constant $c>0$ independent of $p\in [1,\infty]$.

By Lemma \ref{L2.1}  with $f_0=1$ and $f_1= f_{t,r},$ where $f_{t,r}$ is the density of $\mu_{t,r}$ with respect to $\mu$ given in \eqref{FT}, we have
 \beq\label{PP} \W_2(\mu_{t,r},\mu)^2\le \int_M \ff{|\nn L^{-1}(f_{t,r}-1)|^2}{\M(1,f_{t,r})}\d\mu.\end{equation}
 Let
 \beq\label{PX1} \Xi_r(t)= t \int_M |\nn L^{-1} (f_{t,r}-1)|^2 \d\mu,\ \ t,r>0.\end{equation}
 In the next two lemmas, we show that
$$   \lim_{t\to\infty}    \bigg| \E^\nu \Xi_r(t)  - \sum_{i=1}^\infty \ff 2 {\ll_i^2\e^{2r\ll_i}}\bigg|=0,\ \ r>0 $$ holds  for $\nu=h_\nu\mu$ with $\|h_\nu\|_\infty<\infty$,
 and   $\M(1,f_{t,r})$ is close to $1$ for  large $t$, so that \eqref{PP} implies  the desired upper bound estimate \eqref{E1'} for $\E^\nu$ replacing $\E^x$.

 \beg{lem}\label{L2.2}  Assume {\bf (A1)}. There exists a constant $c>0$ such that
   \beq\label{CC}
     \bigg| \E^\nu  \Xi_r(t)- \sum_{i=1}^\infty \ff 2 {\ll_i^2\e^{2r\ll_i}}\bigg|\\
     \le \ff {c\|h_\nu\|_\infty } t \sum_{i=1}^\infty
     \ff 1 {\ll_i^2\e^{2r\ll_i}},\ \  t \ge 1,r>0   \end{equation}    holds for any probability measure $\nu=h_\nu\mu$, and
     \beq\label{CC'}
      \sup_{x\in M} \bigg| \E^x  \Xi_r(t) - \sum_{i=1}^\infty \ff 2 {\ll_i^2\e^{2r\ll_i}}\bigg|
     \le \ff {c  \|P_{r/2}\|_{2\to\infty}^2} t \sum_{i=1}^\infty
     \ff 1 {\ll_i^2\e^{r\ll_i}},\ \  t \ge 1,r>0. \end{equation} \end{lem}

 \beg{proof}   By \eqref{FT} and \eqref{HTT}, we have $\mu(f_{t,r}-1)=0$ and
  $\|f_{t,r}\|_\infty\le \|P_r\|_{1\to\infty}<\infty$.  Consequently,  \eqref{LLL} implies  $(-L)^{-1}(f_{t,r}-1)\in \D(L)$.
 Then the integration by parts formula and the symmetry of $P_s$ in $L^2(\mu)$ yield
  \beq\label{S1}\beg{split}  & \int_M  |\nn L^{-1}(f_{t,r}-1)|^2 \d\mu=- \E^\nu\int_M \big\{L^{-1}(f_{t,r}-1)\big\}\cdot L\big\{ L^{-1}(f_{t,r}-1)\big\}\d\mu \\
  &=    \int_M \big\{(-L)^{-1}(f_{t,r}-1)\big\} (f_{t,r}-1)\d\mu =   \int_0^\infty \d s \int_M (f_{t,r}-1)P_s(f_{t,r}-1)\d\mu \\
  &=  \int_0^\infty\d s \int_M   \big|P_{\ff s 2} f_{t,r}-1 \big|^2\d\mu.\end{split}  \end{equation}
  By \eqref{NHK},  $P_s\phi_i=\e^{-\ll_i s}\phi_i$ and $\mu(\phi_i\phi_j)= 1_{i=j}$, we obtain
  \beq\label{NNK} t\int_M |P_{\ff s 2} f_{t,r}-1|^2\d\mu = \sum_{i=1}^\infty \e^{-\ll_i(2r+s)} |\psi_i(t)|^2.\end{equation}
  Combining this with \eqref{S1} we get
  \beq\label{SX1} \Xi_r(t)= \sum_{i=1}^\infty \ff{|\psi_i(t)|^2}{\ll_i \e^{2\ll_i r}},\ \ t,r>0.\end{equation}
  Moreover,     the Markov property and $P_s\phi_i=\e^{-\ll_i s}\phi_i$ imply
  $$\E^\nu(\phi_i(X_{s_2})|\F_{s_1})= P_{s_2-s_1}\phi_i(X_{s_1}) =\e^{-\ll_i (s_2-s_1)} \phi_i(X_{s_1}),\ \ s_2\ge s_1\ge 0,$$
so that
\beq\label{SX2} \beg{split}  & \E^\nu |\psi_i(t)|^2 =\ff 1 {t} \E^\nu \bigg|\int_0^t  \phi_i(X_s) \d s \bigg|^2
 =   \ff 2 {t} \int_0^t\d s_1 \int_{s_1}^t \E^\nu\big[\phi_i(X_{s_1}) \phi_i(X_{s_2}) \big]\d s_2\\
&=   \ff 2 {t} \int_0^t \E^\nu |\phi_i(X_{s_1})|^2 \d s_1 \int_{s_1}^t   \e^{-\ll_i(s_2-s_1)}\d s_2= \ff 2 {\ll_i t} \int_0^t \nu(P_s\phi_i^2) (1-\e^{-\ll_i(t-s)}) \d s.\end{split}\end{equation}
Combining \eqref{SX1} and \eqref{SX2} gives
 \beq\label{S3}   \beg{split}  &\E^\nu \Xi_r(t) =    \ff {2} {t} \sum_{i=1}^\infty \ff{\e^{-2r \ll_i}}{\ll_i^2} \int_0^t \nu(P_{s} \phi_i^2)  \big(1-\e^{-\ll_i(t-s)}\big)\d s=I_1+I_2,\end{split} \end{equation}  where
\beq\label{S4}   I_1   :=  \ff 2t \sum_{i=1}^\infty \int_0^t \ff { 1-\e^{-(t-s_1)\ll_i}}  {\ll_i^2\e^{2r\ll_i}}  \d s_1
  =  \sum_{i=1}^\infty  \ff {2} {\ll_i^2\e^{2r\ll_i}} -\ff 2t  \sum_{i=1}^\infty  \ff {1-\e^{-\ll_i t}}  {\ll_i^3\e^{2r\ll_i}},   \end{equation} and due to $\nu(P_{s}\phi_i^2)= \mu(h_\nu P_{s}\phi_i^2)=\mu(\phi_i^2P_{s}h_\nu)$,
\beq\label{*PQ} I_2:= \ff 2t \sum_{i=1}^\infty \int_0^t \ff {1-\e^{-(t-s)\ll_i}}  {\ll_i^2\e^{2r\ll_i} }  \mu(\phi_i^2 P_{s}h_\nu -1)   \d s.\end{equation}
Since $\mu(\phi_i^2)=1$, by \eqref{ULT-2} we find a constant $c_1>0$ such that
\beg{align*} |\mu(\phi_i^2P_{s}h_\nu-1)|= |\mu((P_{s }h_\nu-1)\phi_i^2)|\le \|P_{s}(h_\nu-1)\|_\infty\le c_1\e^{-\ll_1 s}\|h_\nu\|_\infty,\ \ s\ge 0.\end{align*}
Thus, there exists a constant $c_2>0$ such that
 $$|I_2|\le    \ff{c_2}t \|h_\nu\|_\infty\sum_{i=1}^\infty\ff 1 {\ll_i^2\e^{2r\ll_i}}<\infty. $$
 Combining this with \eqref{S3} and \eqref{S4}, and noting that $t \ge 1$ and $\|h_\nu\|_\infty\ge 1$, we  prove  \eqref{CC}  for some constant $c>0$.

 Next, when $\nu=\dd_x$ \eqref{S3} becomes
 \beq\label{SS1} \E^\nu \Xi_r(t) \le I_1+I_2(x),\end{equation} where $I_1$ is in \eqref{S4}, and due to $\mu(\phi_i^2)=1$ and $ P_{r/2}\phi_i=\e^{-r\ll_i/2 }\phi_i$,
  \beg{align*} I_2(x)&:= \ff 2t \sum_{i=1}^\infty \int_0^t \ff {1-\e^{-(t-s)\ll_i}}  {\ll_i^2\e^{2r\ll_i} }    P_{s}\big\{\phi_i^2(x)-1 \big\}   \d s\\
&\le \ff 2t \sum_{i=1}^\infty \int_0^t \ff {1}  {\ll_i^2\e^{r\ll_i} }   \big|P_{s}(P_{r/2} \phi_i)^2(x)-\mu((P_{r/2}\phi_i)^2)\big|   \d s.\end{align*}
By \eqref{ULT-2} and noting that $\|P_{s}\phi_i\|_{\infty} \le \|P_s\|_{2\to\infty}$, we find a constant $c_3>0$ such that
\beg{align*}\sup_{x\in M}  I_2(x) &\le    \ff {c_3}t \sum_{i=1}^\infty \int_0^t \ff {1}  {\ll_i^2\e^{r\ll_i} }   \|(P_{r/2} \phi_i)^2\|_\infty \e^{-\ll_1s}   \d s\\
&\le  \ff {c_3\|P_{r/2}\|_{2\to\infty}^2}t \sum_{i=1}^\infty \ff {1}  {\ll_i^2\e^{r\ll_i} }   \int_0^t \e^{- \ll_1s}   \d s\\
&\le \ff {c_3\|P_{r/2}\|_{2\to\infty}^2} {\ll_1 t} \sum_{i=1}^\infty   \ff {1} {\ll_i^2\e^{r\ll_i} }.\end{align*}
     Combining this with \eqref{SS1} and \eqref{S4}, we prove \eqref{CC'} for some constant $c>0$.    \end{proof}

The following lemma is similar to    \cite[Proposition 2.6]{SHI}, which ensures that  $\M(1,f_{t,r})\to 1$ as  $t\to\infty$.

\beg{lem}\label{L2.3}    Assume {\bf (A1)}. Let $\|f_{t,r}-1\|_\infty= \sup_{y\in M}|f_{t,r}(y)-1|$.  Then there exists a function $c: \mathbb N\times (0,\infty)\to (0,\infty)$ such that
$$  \sup_{x\in M} \E^x \big[\|f_{t,r}-1\|_\infty^{2k}\big] \le c(k,r)    t^{-k},\ \ t\ge 1,r>0, k\in \mathbb N. $$
 \end{lem}

\beg{proof} For fixed $r>0$ and $y\in M$, let  $f= p_r(\cdot,y)-1.$ For any $k\in \mathbb N$, consider
$$ I_k(s):= \E^\nu \bigg|\int_0^s f(X_t)\d t\bigg|^{2k}=(2k)!\E^\nu\int_{\DD_k(s)} f(X_{s_1})\cdots f(X_{s_{2k}}) \d s_1\cdots \d s_{2k},\ \ s>0,$$
where $\DD_k(s):= \big\{(s_1,\cdots, s_{2k})\in [0,s]: \ 0\le s_1\le s_2\le \cdots\le s_{2k}\le s\big\}.$
By the Markov property, we have
$$\E^\nu\big(f(X_{s_{2k}})\big|X_t, t\le s_{2k-1}\big) = (P_{s_{2k}-s_{2k-1}}f)(X_{s_{2k-1}}).$$ So, letting $g(r_1,r_2)= (fP_{r_2-r_1}f)(X_{r_1})$ for $r_2\ge r_1\ge 0$,  we obtain
$$ I_k(s) =   (2k)!   \E^\nu \bigg[\int_0^s f(X_{s_1})\d s_1\int_{s_1}^s f(X_{s_2})\d s_2\cdots \int_{s_{2k-2}}^s   \d s_{2k-1}\int_{s_{2k-1}}^s  g(s_{2k-1},s_{2k}) \d s_{2k}\bigg].$$
By the Fubini formula,  we may rewrite $I_k(s)$ as
 \beg{align*} &I_k(s) = (2k)!\, \E^\nu\bigg[\int_{\DD_1(s)} g(r_1,r_2)\d r_1\d r_2 \int_{\DD_{k-1}(r_1)} f(X_{s_1})\cdots f(X_{s_{2k-2}}) \d s_1\cdots \d s_{2k-2}\bigg]\\
&= \ff{(2k)!}{(2k-2)!}  \int_{\DD_1(s)} \E^\nu\bigg[g(r_1,r_2)\bigg|\int_0^{r_1} f(X_r)\d r\bigg|^{2k-2} \bigg]\d r_1\d r_2.\end{align*}
Using  H\"older's inequality, we derive
\beg{align*} &I_k(s) \le 2k(2k-1)  \int_{\DD_1(s)} \big(\E^\nu|g(r_1,r_2)|^k\big)^{\ff 1 k} \bigg(\E^\nu \bigg|\int_0^{r_1} f(X_r)\d r\bigg|^{2k}\bigg)^{\ff{k-1}k} \d r_1\d r_2\\
&\le 2k(2k-1) \bigg(\sup_{u\in [0,s]} I_k(u)\bigg)^{\ff{k-1}k}  \int_{\DD_1(s)} \big(\E^\nu|g(r_1,r_2)|^k\big)^{\ff 1 k} \d r_1\d r_2.\end{align*}
Thus,
$$\sup_{s\in [0,t]} I_k(s) \le 2k(2k-1)    \bigg(\sup_{s\in [0,t]} I_k(s)\bigg)^{\ff{k-1}k}  \int_{\DD_1(t)} \big(\E^\nu|g(r_1,r_2)|^k\big)^{\ff 1 k} \d r_1\d r_2,\ \ t>0.$$
Since $I_k(t)\le (\|f\|_\infty t)^{2k}<\infty$, this implies
\beq\label{**1} I_k(t)\le \sup_{s\in [0,t]} I_k(s)\le \{2k(2k-1)\}^k \bigg(\int_{\DD_1(t)} \big(\E^\nu|g(r_1,r_2)|^k\big)^{\ff 1 k} \d r_1\d r_2\bigg)^{ k}. \end{equation}
Recalling that $g(r_1,r_2)= (fP_{r_2-r_1}f)(X_{r_1})$ and
$$\|f\|_\infty =\|p_r(\cdot,y)-1\|_\infty\le 2\|P_r\|_{1\to\infty}<\infty,$$  by \eqref{ULT-2} we obtain
   $$  |g(r_1,r_2)|^k \le \|f P_{r_2-r_1} f\|_\infty^k \le c\e^{-\ll_1(r_2-r_1)k} \|f\|_\infty^{2k} \le c \|P_r\|_{1\to\infty}^{2k} \e^{-\ll_1(r_2-r_1)k}$$ for some constant $c>0$.
  Thus,
 \beg{align*}  &  \bigg(\int_{\DD_1(t)} \big(\E^x|g(r_1,r_2)|^k\big)^{\ff 1 k} \d r_1\d r_2\bigg)^{ k}\\
&\le
\bigg(\int_0^t   \d r_1 \int_{r_1}^t c\|P_r\|_{1\to\infty}^2\e^{- \ll_1(r_2-r_1)}\d r_2\bigg)^k\\
&\le \big(c\ll_1^{-1}\|P_r\|_{1\to\infty}^2 t\big)^k,  \ \ t\ge 1, r>0,  k\in\mathbb N.
\end{align*}  This and \eqref{**1} yield
\beq\label{L2.310} \sup_{x,y\in M} \E^x\big[ |f_{t,r}(y)-1|^{2k}\big]= t^{-2k} I_k(t) \le c(k,r) \|P_r\|_{1\to\infty}^{2k}   t^{-k},\ \ t\ge 1,r>0 \end{equation} for all $k\in \mathbb N$ and some constant  $c(k)>0$.

Finally,   noting that $f_{t,r}=P_{r/2} f_{t,r/2},$ we deduce from \eqref{L2.310} that
\beg{align*} &\sup_{x\in M} \E^x \big[\|f_{t,r}-1\|_\infty^{2k}\big] = \sup_{x\in M} \E^x \big[\|P_{\ff r 2} (f_{t,\ff r 2}-1)\|_\infty^{2k}\big] \\
&\le \|P_{\ff r 2}\|_{2k\to \infty}^{2k}\sup_{x\in M} \E^x\big[\mu(|f_{t,\ff r 2}-1|^{2k})\big]
 \le c(k) \|P_{\ff r 2}\|_{1\to\infty}^{4k} t^{-k},\ \ t\ge 1, r>0.\end{align*} This finishes the proof.
\end{proof}

We are now ready to prove the upper bound estimate  \eqref{E1'} in Theorem \ref{T1.2}.

\beg{prp}\label{P1}  The assumption {\bf (A1)} implies   $\eqref{E1'}$.  \end{prp}

\beg{proof} (a) Proof of \eqref{E1'}. By \eqref{HK}, \eqref{HTT} and $\mu(\phi_i^2)=1$, we have
 $$\sum_{i=1}^\infty\ff 1 {\ll_i^2\e^{2r\ll_i}} \le \ff 1 {\ll_1^2} \sum_{i=1}^\infty \e^{-2r\ll_i}=\ff 1 {\ll_1^2} \int_M \big( p_{2r}(x,x) - 1 \big) \mu(\d x)\le \ff{\|P_{2r}\|_{1\to\infty}}{\ll_1^2}<\infty.$$ So, it remains to prove the first inequality in \eqref{E1'}.

For any  $\eta\in (0,1)$, consider  the event
\beq\label{AN} A_\eta= \Big\{  \|f_{t,r}-1\|_\infty \le \eta\Big\}.\end{equation}
Noting that   $f_{t,r}(y)\ge 1-\eta$ implies
$$\M(1, f_{t,r}(y))\ge \ss{f_{t,r}(y)}\ge \ss{1-\eta},$$
we deduce from   Lemma \ref{L2.1} and \eqref{CC'} that  for some constant $c(r)>0$,
\beg{align*}
&t\sup_{x\in M} \E^x [1_{A_\eta}\W_2(\mu_{t,r},\mu)^2]\le \sup_{x\in M} \E^x\bigg\{\ff{\Xi_r(t)}{  \ss{1-\eta} }\bigg\} \\
&\le  \ff 1{\ss{1-\eta}} \sum_{i=1}^\infty \ff{2}{\ll_i^2\e^{2\ll_i r}} \Big(1 + \ff {c(r)}t \Big),\ \ t>0, \eta\in (0,1).\end{align*}
So,
\beq\label{N1}\beg{split} &   t\sup_{x\in M}\E^x [\W_2(\mu_{t,r},\mu)^2]\le  \ff 1{\ss{1-\eta}} \sum_{i=1}^\infty \ff{2}{\ll_i^2\e^{2\ll_i r}} \Big(1 + \ff {c(r)}t \Big)+t\sup_{x\in M} \E^x\big[1_{A_\eta^c} \W_2(\mu_{t,r},\mu)^2\big]  \\
&\le  \ff {1+c(r)t^{-1}}{\ss{1-\eta}} \sum_{i=1}^\infty \ff{2}{\ll_i^2\e^{2\ll_i r}} + t \sup_{x\in M} \ss{ \P^x (A_\eta^c)\E^x [\W_2(\mu_{t,r},\mu)^4]},\ \ t, \eta\in (0,1).\end{split}\end{equation}
By Jensen's inequality and \eqref{MM},  we obtain
\beq\label{NN0} \beg{split} &\E^x \W_2(\mu_{t,r},\mu)^4\le \E^x \bigg(\int_{M\times M} \rr(z,y)^2\mu_{t,r}(\d z)\mu(\d y)\bigg)^2\\
& \le \E^x \int_{M\times M} \rr(z,y)^4\mu_{t,r}(\d z)\mu(\d y)\le \ff 1 t\int_0^t \E^x \mu\big(\rr(X_{r+s},\cdot)^4 \big)\d s\\
&\le  \ff 1 t\int_0^t  \|P_{s+r}\|_{1\to\infty} (\mu\times \mu) (\rr^4)\d s\le \|P_r\|_{1\to\infty}   (\mu\times\mu)(\rr^4)<\infty.\end{split}\end{equation}
Moreover, Lemma \ref{L2.3} implies
 \beq\label{*D} \sup_{x\in M}\P^x (A_\eta^c) \le \eta^{-2k} c(k,r) t^{-k},\ \   t\ge 1, k\in \mathbb N, \eta\in (0,1) \end{equation} for some constant $c(k,r)>0$.
By taking $k=4$ in \eqref{*D} and applying  \eqref{N1} and \eqref{NN0},   we  conclude that
$$ \limsup_{t\to\infty} \Big\{t\sup_{x\in M} \E^x   [\W_2(\mu_{t,r},\mu)^2]\big\} \le
  \ff 1 {\ss{1-\eta}}\sum_{i=1}^\infty\ff{2}{\ll_i^2\e^{2r\ll_i}}, \ \ \eta\in (0,1). $$ By letting $\eta\downarrow 0$, we  derive \eqref{E1'}.
\end{proof}

\subsection{Lower bound estimate}

  Due to \eqref{E1'}, \eqref{E2'} follows from the lower bound estimate
 \beq\label{E2-3} \liminf_{t\to\infty} \Big\{t\inf_{x\in M} \E^x \W_2(\mu_{t,r},\mu )^2\Big\} \ge \sum_{i=1}^\infty \ff 2 {\ll_i^2\e^{2r\ll_i}},\ \ r>0.\end{equation}
To estimate $\W_2(\mu_{t,r},\mu)$ from below, we use the fact that
 \beq\label{LBD} \beg{split} &\ff 1 2 \W_2(\mu_{t,r},\mu)^2\ge \mu_{t,r}(\phi_1)- \mu(\phi_0),\ \ (\phi_0,\phi_1)\in \C,\\
 &\C:=\Big\{(\phi_0,\phi_1): \phi_0,\phi_1\in  C_b(M),   \phi_1(x)-\phi_0(y)\le \ff 12\rr(x,y)^2 \ \text{for}\  x,y\in M\Big\}.\end{split}\end{equation}
  We will construct the pair $(\phi_0,\phi_1)$ by using the idea of \cite{AMB}, where compact $M$ without boundary has been considered.  To realize the idea in the present more general setting,
  we need the following result on gradient estimate which is implied by \cite[Corollary 1.2]{W14a} for $Z=\nn V$.

  \beg{lem}[\cite{W14a}] \label{L009} If there exists $\phi\in C_b^2(M)$ such that $\inf \phi=1, |\nn\phi | \cdot |\nn V|$ is bounded,   $\nn \phi\parallel N ($i.e. $\nn\phi$ is parallel to $N)$ and $ \II\ge -N\log \phi$ hold on $\pp M$, and
  $$\Ric_V-\ff 1 2 \phi^2  L\phi^{-2} \ge -K_\phi $$ holds for some constant $K_\phi\ge 0$.
  Then
\beg{align*} & |\nn P_t f|^2\le \ff {\e^{2K_\phi t}}{\phi^2} P_t(\phi|\nn f|)^2,\ \ t\ge 0, f\in C_b^1(M),\\
 & |\nn P_t f|^2\le \ff {\|\phi\|_\infty K_\phi }{\e^{2K_\phi t}-1} \big\{P_t(f^2)-(P_tf)^2\big\},\ \ t>0, f\in \B_b(M),\\
 &  P_t (f^2) \le (P_t f)^2 + \ff{\|\phi\|_\infty^2(\e^{2K_\phi t}-1)}{K_\phi}P_t|\nn f|^2,\ \ t>0, f\in C_b^1(M).\end{align*}
 \end{lem}
As a consequence of Lemma \ref{L009}, we have the following result.

 \beg{lem}\label{L3.0} Assume {\bf (A2)}. There exists a constant $\kk>0$ such that
 \beq\label{W141} |\nn P_t f|^2\le \big(1+\kk\ss t\big) P_t|\nn f|^2,\ \ t\in [0,1], f\in C_b^1(M),\end{equation}
\beq\label{W142} |\nn P_t f|^2\le \ff \kk{t} P_t(f^2),\ \ t\in (0,1], f\in \B_b(M),\end{equation}
\beq\label{W14P}P_t (f^2) \le (P_t f)^2 + \kk t P_t|\nn f|^2,\ \ t\in (0,1], f\in C_b^1(M).\end{equation}
  Consequently,
 \beq\label{ERR} P_t \e^{4 f}   \le 8 (P_t \e^{f})^4,\ \ t \in (0,1], \|\nn f\|_\infty^2\le \ff 1 {8\kk t}.\end{equation} \end{lem}

 \beg{proof} Let $\Ric_V\ge -K$ for some constant $K\ge 0$. If $\pp M$ is empty or convex, we may take $\phi=1$ and $K_\phi=K$ in Lemma \ref{L009}. Then
\eqref{W141}-\eqref{W14P} follow from     estimates  in Lemma \ref{L009}.

 If $\pp M\ne\emptyset$ and there exists a compact set $D$ such that $\pp M$ is convex outside $D$, we make use of Lemma \ref{L009}.
 To this end, we construct a function $g\in C^\infty_0(M)$ such that $0\le g\le 1$,   $Ng|_{\pp M}=0$, and $g=1$ on the compact set $D$. Let $D'$ be the support of $g$.
 Since the distance $\rr_\pp$ to the boundary is smooth in a neighborhood of $\pp M$, we may take a constant $r_0\in (0,1)$  such that  $\rr_\pp$ is smooth on $D'\cap \pp_{r_0}M$, where $\pp_{r_0}M:=\{\rr_\pp\le r_0\}\subset M$. Moreover, since $\I$ is nonnegative on $\pp M\setminus D$, there exists a constant $\kk>0$ such that $\I\ge -\kk$. We choose $h\in C^\infty([0,\infty))$ such that $h$ is increasing,
 $h(r)=r$ for $r\in [0,\ff {r_0} 2]$ and $h(r)= h(r_0)$ for $r\ge r_0$. For any $\vv\in (0,1)$, take
 $$\phi = 1+ \kk \vv g h(\vv^{-1}\rr_\pp).$$
 It is easy to see that $\inf \phi=1,  \nn \phi\parallel N$ and $ \II\ge -N\log \phi$ hold on $\pp M$ as  required by Lemma \ref{L009}. Next, since $\phi\ge 1$ and $\nn\phi=0$ outside the compact set $D'$, there exists a constant $c_1>0$ such that
 $$\ff 1 2 \sup_M\{\phi^2 L\phi^{-2}\}= \sup_{D'} \{3\phi^{-2}|\nn \phi|^2 -\phi^{-1} L\phi\} \le c_1\vv^{-1},\ \ \vv\in (0,1).$$
 Combining this with \eqref{CV}, we obtain
 $$\Ric_V- \ff 1 2 \phi^2 L\phi^{-2}\ge - K - c_1\vv^{-1}\ge -c_2\vv^{-1},\ \ \vv\in (0,1)$$ for some constant $c_2>0$.
 Then the second and third estimates follow  from \eqref{W142} and \eqref{W14P}, while \eqref{W141} implies
 \beg{align*}& |\nn P_t f|^2\le \ff{\e^{2c_2\vv^{-1}t} }{\phi^2} P_t(\phi|\nn f|)^2\le \e^{2c_2\vv^{-1}t}  \|\phi\|_\infty^2 P_t |\nn f|^2\\
 &\le
 \e^{2c_2\vv^{-1}t}  (1+ \kk\|h\|_\infty\vv)^2 P_t |\nn f|^2,\ \ t,\vv\in (0,1).\end{align*}
 Taking $\vv= \ss t,$ we prove the first estimate for some constant $c>0$.

 It remains to prove \eqref{ERR}. By \eqref{W14P}, we have
 $$P_t \e^{2f}  \le (P_t \e^{f})^2 + \kk t P_t(|\nn f|^2 \e^{2f})\le  (P_t \e^{f})^2 + \kk t \|\nn f\|_\infty^2 P_t( \e^{2f}).$$
 This implies
 $$P_t\e^{2f}\le 2 (P_t \e^{f})^2,\ \ \text{if}\ \|\nn f\|_\infty^2\le \ff 1{2\kk t}.$$
Using $2f$ replacing $f$ we obtain
  $$P_t\e^{4f}\le 2 (P_t \e^{2f})^2,\ \ \text{if}\ \|\nn f\|_\infty^2\le \ff 1{8\kk t}.$$
  Therefore, \eqref{ERR} holds.
 \end{proof}

We are now ready to present the following key lemma for the lower bound estimate of $\W_2(\mu_{t,r},\mu)$.

\beg{lem}\label{L3.1} Assume {\bf (A1)} and  {\bf (A2)}.
 For any $f\in C^2_b(M)$ with $\|\nn f\|_\infty+\|Lf\|_\infty<\infty$ and $Nf|_{\pp M}=0$ if $\pp M\ne \emptyset$,  let
$\phi_t^\si= -\si\log P_{\ff{\si t} 2} \e^{-\si^{-1}f },\ \ t\in [0,1],\ \si>0.$ Then $\phi_t^\si\in C^2(M)$ and
\beg{enumerate} \item[$(1)$] $\phi_0^\si=f, \|\phi_t^\si\|_\infty\le \|f\|_\infty,$ and $\pp_t \phi_t^\si= \ff\si 2 L \phi_t^\si -\ff 1 2 |\nn \phi_t^\si|^2, t>0;$
\item[$(2)$]   There exist constants $c,\gg>0$ such that for any $\si,t\in (0,1],$ when $\|\nn f\|_\infty^2\le \gg\si$ we have
    \beg{align*} &\phi_1^\si(y)-\phi_0^\si(x)\le \ff 1 2\Big\{\rr(x,y)^2 +     \si\|(Lf)^+\|_\infty+ c  \si^{\ff 1 2}  \|\nn f\|_\infty^2\Big\},\\
&\int_M (\phi_0^\si-\phi_1^\si) \d\mu \le \ff 1 2       \int_M  |\nn f|^2 \d\mu + c\si^{-1}\|\nn f\|_\infty^4.\end{align*}
\end{enumerate}
\end{lem}

\beg{proof}   (1)  The first assertion follows from standard calculations. Indeed,  by the chain rule and the heat equation   $\pp_t P_t g= LP_tg$ for $t>0$ and $g\in C_b(M)$, we have
$$\pp_t\phi_t^\si = -\ff{\si^2LP_{\ff{t\si}2}\e^{-\si^{-1} f}}{2P_{\ff{t\si}2}\e^{- \si^{-1} f}} =\ff \si 2 L \phi_t^\si -\ff 1 2 |\nn \phi_t^\si|^2.$$

(2)   Let $\si,t\in (0,1]$ and $\|\nn f\|_\infty^2\le \gg \si$ for $\gg= \ff 1{4\kk}$, where $\kk>0$ is in   Lemma \ref{L3.0}. Then
$\|\si^{-1}\nn f\|_\infty^2\le \ff{1}{4\kk\si}\le \ff 1 {8 \kk t'}$ for $t'=\ff{t\si}2$, so that
 \eqref{ERR} holds for $(t',-\si^{-1}f)$ replacing $(t,f)$. This  implies
\beq\label{**00}  P_{\ff {t\si}2} \e^{-4\si^{-1}f} \le 8 (P_{\ff {t\si}2} \e^{-\si^{-1}f})^4.\end{equation} Combining   with \eqref{W141} gives
\beq\label{*ZZ} |\nn\phi_t^\si|^2 =\ff{|\nn P_{\ff{t\si}2} \e^{-\si^{-1}f}|^2}{(P_{\ff{t\si}2} \e^{-\si^{-1}f})^2} \le \ff{(1+\kk)   P_{\ff{t\si}2} (|\nn f|^2 \e^{-2\si^{-1}f})}{(P_{\ff{t\si}2} \e^{-\si^{-1}f})^2}  \le c \|\nn f\|_\infty^2\end{equation} for some constant $c>0$.

Next, by (2.1) in \cite{W05}, for $g \in C_b^1(M)$,
\beq\label{gra} |\nn P_t g|(x) \le \E^x \big[ | \nn g|(X_t) \e^{Kt + \dd l_t } \big] ,\ \ x \in M,\ \ t>0, \end{equation}
where $K\ge 0$ is the constant such that $\Ric_V\ge -K$  holds on $M$, $\dd > 0$ is the constant such that $\I\ge -\dd$ and $l_t$ is the local time of the $L$-process on the boundary $\pp M$.
Combining this with \eqref{**00},
 $$LP_{\ff{t\si}2} \e^{-\si^{-1}f} = P_{\ff{t\si}2} L\e^{-\si^{-1}f} =-\ff 1 \si P_{\ff{t\si}2} (\e^{-\si^{-1}f} Lf) +\ff 1 {\si^2} P_{\ff{t\si}2} (|\nn f|^2\e^{-\si^{-1}f}),$$ and  applying H\"older's inequality, we find a constant $c_0>0$ such that
 \beq\label{ASC'} \beg{split}&L\phi_t^\si=-\ff{\si LP_{\ff{t\si}2} \e^{-\si^{-1}f}}{P_{\ff{t\si}2}\e^{-\si^{-1}f}} +\ff{\si|\nn P_{\ff{t\si}2}\e^{-\si^{-1}f}|^2}{(P_{\ff{t\si}2}\e^{-\si^{-1}f})^2} \\
 &\le \|(Lf)^+\|_\infty-\ff{P_{\ff{t\si}2}(|\nn f|^2\e^{-\si^{-1}f})}{\si P_{\ff{t\si}2}\e^{-\si^{-1}f}} + \ff{\big( \E \big[(\e^{-\si^{-1}f}| \nn f|)(X_{\ff{t\si}2}) \e^{\ff{Kt\si}2 + \dd l_{\ff{t\si}2} } \big]\big)^2}{\si (P_{\ff{t\si}2}\e^{-\si^{-1}f})^2}\\
 &\le \|(Lf)^+\|_\infty-\ff{P_{\ff{t\si}2}(|\nn f|^2\e^{-\si^{-1}f})}{\si P_{\ff{t\si}2}\e^{-\si^{-1}f}} + \ff{P_{\ff{t\si}2}(|\nn f|^2\e^{-\si^{-1}f}) \E\big[\e^{-\si^{-1}f(X_{\ff{t\si}2}) + Kt\si + 2\dd l_{\ff{t\si}2} }\big]}{\si (P_{\ff{t\si}2}\e^{-\si^{-1}f})^2}\\
 & = \|(Lf)^+\|_\infty + \ff{P_{\ff{t\si}2}(|\nn f|^2\e^{-\si^{-1}f}) \E \big[\e^{-\si^{-1}f(X_{\ff{t\si}2})}   \big(\e^{Kt\si + 2\dd l_{\ff{t\si}2} }-1\big)\big]}{\si (P_{\ff{t\si}2}\e^{-\si^{-1}f})^2}\\
 &\le \|(Lf)^+\|_\infty +  c_0 \si^{-1} \|\nn f\|_\infty^2  \Big(  \E \big[ (\e^{Kt\si + 2\dd l_{\ff{t\si}2} }-1)^{\ff 4 3}\big]\Big)^{\ff 3 4}.\end{split} \end{equation}
By the proof of Lemma 2.1 in \cite{W05},  for any $\ll > 0$  there exists a constant $c>0$ such that
\begin{align*}
\E \e^{\ll l_t} \le c(\ll), \, \ t\in (0, 1].
\end{align*}
Moreover,  by Lemma 2.2 in \cite{W09'},   there exists a constant $c_1>0$ such that
\begin{align*}
\E l_t^2 \le c_1t, \, \ t\in (0, 1].
\end{align*}
Combining these facts, we find a constant $c= c(K, \dd)>0$ such that
\beq\label{lt} \beg{split}&\E \big[ (\e^{Kt\si + 2 \dd l_{\ff{t\si}2} } - 1)^{\ff 4 3}\big] \le   \E \big[ (Kt\si+2\dd l_{\ff{t\si}2})^{\ff 4 3} \e^{\ff 4 3 Kt\si + \ff 8 3 \dd l_{\ff{t\si}2} } \big]  \\
& \le \Big(\E [(Kt\si+2\dd l_{\ff{t\si}2})^2]\Big)^{\ff 2 3}\Big(\E \e^{ 4  Kt\si +   8  \dd l_{\ff{t\si}2} }\Big)^{\ff 1 3} \le    c    \si^{\ff 2 3},\  \ \si,t\in (0,1].\end{split} \end{equation}
 It then follows from  \eqref{ASC'} and \eqref{lt} that
\beq\label{ASC}  L\phi_t^\si \le \|(Lf)^+\|_\infty + c_2  \si^{-\ff12}  \|\nn f\|_\infty^2,\ \ \si,t\in (0,1], \|\nn f\|_\infty^2\le \gg\si  \end{equation} holds for some constant $c_2>0$.

Now, for any two points $x,y\in M$,  let $\gg: [0,1]\to M$ be the minimal geodesic from $x$ to $y$,
 so that $|\dot\gg_t|=\rr(x,y)$. By (1) and  \eqref{ASC}, we derive
  \beg{align*}  &\ff{\d }{\d t} \phi_t^\si (\gg_t) = (\pp_t \phi_t^\si) (\gg_t) + \<\nn \phi_t^\si (\gg_t), \dot\gg_t\> \\
 &= -\ff 1 2 |\nn \phi_t^\si(\gg_t)|^2 +\ff \si 2 L\phi_t^\si(\gg_t) +\<\nn \phi_t^\si(\gg_t), \dot\gg_t\> \\
 &\le \ff 1 2 |\dot\gg_t|^2 +\ff\si 2 \|(Lf)^+\|_\infty +  c    \ss\si  \|\nn f\|_\infty^2 \\
 &=\ff 1 2 \rr(x,y)^2 +\ff\si 2 \|(Lf)^+\|_\infty + c  \ss\si     \|\nn f\|_\infty^2,\ \ \si,t\in [0,1], \|\nn f\|_\infty^2\le \gg\si  \end{align*} for some constant $c>0$.
 Integrating over $t\in [0,1]$ and noting that $\phi_0^\si(x)=f(x)$, we derive   the first inequality in (2).

 On the other hand,  since $\phi_t^\si \in C^2(M)$ with $N\phi_t^\si|_{\pp M}=0$ and bounded $|\nn\phi_t^\si|+|L\phi_t^\si|$, we have $\mu(L\phi_t^\si)=0$ so that assertion (1)   yields
 \beq\label{PD1}\beg{split}& \mu(f-\phi_1^\si)= \int_M(\phi_0^\si-\phi_1^\si)\d\mu =-\int_M\d\mu\int_0^1(\pp_t \phi_t^\si)\d t\\
 &=\int_0^t \d t \int_M\Big\{\ff 1 2 |\nn \phi_t^\si|^2 -\ff\si 2 L\phi_t^\si\Big\}\d\mu =\ff 1 2 \int_0^1 \mu(|\nn \phi_t^\si|^2) \d t.\end{split} \end{equation}
 Since $\phi^\si \in C^2((0,\infty)\times M)$ with $N\phi_s^\si|_{\pp M}=0$ for $s>0$, we have
 $$N\pp_s\phi_s^\si|_{\pp M}= \pp_s N\phi_s^\si|_{\pp M}=0.$$ Combining this with assertion (1) and applying
 the integration by parts formula, we obtain
 \beg{align*}&\ff{\d}{\d s} \mu(|\nn\phi_s^\si|^2) =-\ff{\d}{\d s} \int_M \phi_s^\si L \phi_s^\si \d\mu =-\int_M (L\phi_s^\si)\pp_s\phi_s^\si\d\mu-\int_M \phi_s^\si L(\pp_s\phi_s^\si)\d\mu\\
 &=-2 \int_M (L\phi_s^\si)\pp_s \phi_s^\si\,\d\mu= - 2 \int_M (L\phi_s^\si)\Big(\ff \si 2 L\phi_s^\si -\ff 1 2 |\nn \phi_s^\si|^2\Big)\d\mu,\ \ s>0.  \end{align*}
  This and \eqref{*ZZ}   imply
 \beq\label{*X3}\beg{split} &\mu(|\nn\phi_t^\si|^2)-\mu(|\nn f|^2)= \int_0^t \Big\{\ff{\d}{\d s} \mu(|\nn\phi_s^\si|^2)\Big\}\d s \\
 &=-2 \int_0^t\d s \int_M (L\phi_s^\si)\Big(\ff \si 2 L\phi_s^\si -\ff 1 2 |\nn \phi_s^\si|^2\Big)\d\mu\\
  &\le   \ff 1 {4\si} \int_0^t \mu(  |\nn \phi_s^\si|^4) \d s \le c \si^{-1}  \|\nn f\|_\infty^4,\ \ \si, t\in [0,1],\|\nn f\|_\infty^2\le \gg\si \end{split}\end{equation}
 for some constant $c>0$.   Substituting this into \eqref{PD1}, we prove the second estimate in assertion (2).

\end{proof}

We are now ready to prove the estimate \eqref{E2-3}.

\beg{prp}\label{P2} Assumptions {\bf (A1)} and {\bf (A2)} imply  $\eqref{E2-3}$.\end{prp}

\beg{proof} Let $f=L^{-1}(f_{t,r}-1)$, and denote
  \beg{align*} &C_1(f,\si):=   \si^{-1} \|\nn f\|_\infty^4,\\
  &C_2(f,\si):= \si\|f_{t,r}-1\|_\infty  + c \si^{\ff 1 2}  \|\nn f\|_\infty^2, \end{align*}
  where $c>0$ is the constant in Lemma \ref{L3.1}(2).
Then
\beq\label{LLO} \|Lf\|_\infty= \|f_{t,r}-1\|_\infty,\end{equation}
and   by \eqref{ULT-2} there exists a constant $c_1>0$ such that
\beq\label{*LLO} \|f\|_\infty \le  \int_0^\infty \|P_s(f_{t,r}-1)\|_\infty\d s \le c_1\|f_{t,r}-1\|_\infty \int_0^\infty \e^{-\ll_1 s}\d s =\ff{c_1}{\ll_1} \|f_{t,r}-1\|_\infty.\end{equation}
Moreover,  by Lemma \ref{L3.0}, there exists a constant $c_0>0$ such that
\beq\label{POO} \|\nn P_tg\|_\infty \le   c_0(1+t^{-\ff 1 2})  \|g\|_\infty,\ \ t>0, g\in \B_b(M). \end{equation}
Combining this with \eqref{ULT-2} implied by
  {\bf (A1)}, we find constants $c_2,c_3,c_4>0$ such that
\beq \label{*DF}\beg{split} &\|\nn f\|_\infty=\|\nn L^{-1} (f_{t,r}-1)\|_\infty \le  \int_0^\infty \|\nn P_s (f_{t,r}-1)\|_\infty\d s,\\
&\le  c_2 \int_0^\infty  (1+s^{-\ff 1 2 }) \|P_{\ff s 2} (f_{t,r}-1)\|_\infty \d s\\
&\le c_3\|f_{t,r}-1\|_\infty \int_0^\infty (1+s^{-\ff 1 2 }) \e^{-\ll_1s/2}\d s\le c_4  \|f_{t,r}-1\|_\infty. \end{split}
\end{equation}
 Let  $B_\si:=\{\|f_{t,r}-1\|_\infty\le \si^{\ff23} \}$. By \eqref{*DF}, there exists a constant $\si_0\in (0,1]$ such that  $B_\si \subset \{\| \nn f\|_\infty^2\le \gg\si \}$ holds for $\si\le \si_0$. So, we deduce
from   \eqref{LLO}, \eqref{*LLO}  and \eqref{*DF}  that
\beq\label{ANN} C_1(f,\si)1_{B_\si} \le c_5\si^{\ff53},\ \ C_2(f,\si)1_{B_\si} \le c_5\si^{\ff53},\ \ \si\in (0,\si_0]  \end{equation}  holds for some constant $c_5>0$.

On the other hand, it is easy to see that $f$ satisfies the Neumann boundary condition, so that by \eqref{LLO} and \eqref{*DF},
Lemma \ref{L3.1}  applies. By Lemma \ref{L3.1}(2),   the integration by parts formula and noting that $f=L^{-1}(f_{t,r}-1)$, we obtain that on the event $B_\si$,
  \beq\label{AN3} \beg{split}  & C_2(f,\si)+\ff 1 2 \W_2 (\mu_{t,r},\mu)^2 \ge \int_M \phi_1^\si\d\mu- \int_M f\d\mu_{t,r}\\
&= \int_M(\phi_1^\si -f)\d\mu- \int_M f(f_{t,r}-1)\d\mu \\
 & \ge -\ff 1 2  \int_M  |\nn L^{-1}(f_{t,r}-1)|^2  \d\mu - \int_M (f_{t,r}-1) L^{-1} (f_{t,r}-1)\d\mu -C_1(f,\si)\\
&= \ff 1 2  \int_M  |\nn L^{-1}(f_{t,r}-1)|^2  \d\mu -C_1(f,\si),\ \ \si\in (0,\si_0].\end{split}\end{equation}
Since $\W_2 (\mu_{t,r},\mu)^2\ge 0$, we deduce from  this,  \eqref{LLO} and \eqref{*DF}  that on the event $B_\si$,
$$\ff 1 2 \W_2 (\mu_{t,r},\mu)^2  \ge    \ff 1 2 \int_M  |\nn L^{-1}(f_{t,r}-1)|^2  \d\mu- C_1(f,\si)-C_2(f,\si).$$
This and \eqref{ANN} yield
\beq\label{*PW2} \beg{split} &\ff t 2\inf_{x\in M} \E^x [ \W_2 (\mu_{t,r},\mu)^2]
\ge \ff t 2 \inf_{x\in M} \E^x [\W_2 (\mu_{t,r},\mu)^21_{B_\si}] \\
\ge &\ff 1 2 \inf_{x\in M} \E^x  \bigg[1_{B_\si}   \Xi_r(t)\bigg]-c_6\si^{\ff53}t\\
&\ge\ff 1 2 \inf_{x\in M} \E^x[\Xi_r(t) ]- I-c_6\si^{\ff53}t,  \ \ \si\in (0,\si_0]\end{split}\end{equation} for some constant $c_6>0$,
where, by \eqref{*DF},   Lemma \ref{L2.3} and noting that $\|f_{t,r}-1\|_\infty\le \|P_r\|_{1\to \infty}<\infty$,
\beq\label{*PW3} \beg{split} &I:= t\sup_{x\in M } \E^x\big[1_{B_\si^c}
\mu\big(  |\nn L^{-1}(f_{t,r}-1)|^2 \big)\big]\\
&\le c_3^2 \|P_r\|_{1\to\infty}^2t \sup_{x\in M} \P^x(B_\si^c)\le  \si^{-\ff{4k}3}  c(k,r)  t^{1-k},\ \ k\in \mathbb N, r>0,\end{split}  \end{equation} where $c(k,r)>0$ is a constant depending on $k,r$.
Now, let $\si=t^{-\aa}$ for some $\aa\in (\ff 35, \ff 34)$ and take $k\ge 1$ such that $k(1- \ff{4\aa}3)>1$. Then we derive from \eqref{*PW2} and \eqref{*PW3} that
\beg{align*} \ff 1 2 \liminf_{t\to\infty}\Big\{  t \inf_{x\in M} \E^x [ \W_2 (\mu_{t,r},\mu)^2]\Big\}
 \ge  \ff 1 2  \liminf_{t\to\infty}  \inf_{x\in M} \E^x \big[ \Xi_r(t)\big].   \end{align*}
Combining this with  \eqref{CC'}, we  prove \eqref{E2-3}. \end{proof}

\subsection{Proofs of Theorem \ref{T1.2} and Corollary \ref{C1.2}}
Since \eqref{E1'} and \eqref{E2'} in Theorem \ref{T1.2} follow from Proposition \ref{P1} and Proposition \ref{P2}, below we only   prove \eqref{CLT} and Corollary \ref{C1.2}.
To this end, we first present the following two lemmas.

\beg{lem}\label{L01} Assume {\bf (A1)}. Then for any $r>0$,  $\Xi_r(\cdot)$ in $\eqref{PX1}$ satisfies
\beq\label{CLT1}  \lim_{t\to \infty}\sup_{\|h_\nu\|_\infty\le C } \big|\P^\nu\big(\Xi_r(t) <a\big)-\nu_r\big((-\infty,a)\big)\big|=0,\ \ a\in \R, C>0.\end{equation}
If $M$ is compact and $d\le 3$, then for any $r_t\downarrow 0$ as $t\uparrow \infty$,
\beq\label{CLT2} \lim_{t\to \infty}\sup_{\|h_\nu\|_\infty\le C} \big|\P^\nu\big(\Xi_{r_t}(t) <a\big)-\nu_0\big((-\infty,a)\big)\big|=0,\ \ a\in \R, C>0.\end{equation}\end{lem}

\beg{proof} By \eqref{SX1} and \eqref{SX2} we have
 \beq\label{GN*0} \Xi_r(t)=\sum_{k=1}^\infty \int_0^\infty \e^{-\ll_k(2r+s)} |\psi_k(t)|^2 \d s = \sum_{k=1}^\infty \ff{|\psi_k(t)|^2}{\ll_k\e^{2\ll_k r}},\ \ t>0.\end{equation}
For any $n\ge 1$, consider the $n$-dimensional process
 $$\Psi_n(t):=(\psi_1(t),\cdots,\psi_n(t)),\ \ t>0.$$ For any $\aa\in \R^n$, we have
 $$\<\Psi_n^{(n)},\aa\>=\ff 1 {\ss t}\int_0^t \Big(\sum_{k=1}^n \aa_k \phi_k(X_s)\Big)\d s.$$
  By \cite[Theorem $2.4'$]{Wu},  when $t\to\infty$, the law of
  $\<\Psi_n(t),\aa\>$ under  $\P^\nu$   converges weakly to the Gaussian distribution $N(0, \si_{n,\aa})$ uniformly in $\nu$ with $\|h_\nu\|_\infty\le C$, where, due to \eqref{SX2} with $\nu=\mu$ and $\mu(P_s\phi_i^2)=\mu(\phi_i^2)=1$,
  the variance is given by
   \beg{align*} &\si_{n,\aa}:= \lim_{t\to\infty} \E^\mu \<\Psi_n(t),\aa\>^2 \\
   &=  \lim_{t\to\infty} \ff 2 t \sum_{k=1}^n \aa_k^2 \int_0^t \d s_1 \int_{s_1}^t \e^{-\ll_k(s_2-s_1)}\d s_2 = \sum_{k=1}^n \ff{2\aa_k^2}{\ll_k}. \end{align*} Thus, uniformly in $\nu$ with $\|h_\nu\|_\infty\le C$,
   $$\lim_{t\to\infty} \E^\nu \e^{{\rm i} \<\Psi_n(t),\aa\>} = \int_{\R^n} \e^{{\rm i} \<x,\aa\>}  \prod_{k=1}^n  N(0,  2\ll_k^{-1})(\d x_k),\ \ \aa\in\R^n,$$ so that
   the distribution of $\Psi_n(t) $ under $\P^\nu$ converges weakly to $\prod_{k=1}^n N(0,  2\ll_k^{-1}).$ Therefore,  letting
  \beq\label{XIN} \Xi_r^{(n)}(t):= \sum_{k=1}^n \ff{|\psi_k(t)|^2}{\ll_k^2\e^{2\ll_k r}},\ \ \Xi_r^{(n)}:=  \sum_{k=1}^n \ff{2\xi_k^2}{\ll_k^2\e^{2\ll_k r}},\end{equation}
we derive
   \beq\label{GN*1} \lim_{t\to\infty} \sup_{\|h_\nu\|_\infty\le C } \big|\P^\nu\big(\Xi_r^{(n)}(t) <a\big)-\P\big(\Xi_r^{(n)} <a\big)\big|=0,\ \ a\in \R.\end{equation}
  On the other hand,  \eqref{SX1}, \eqref{SX2} and \eqref{XIN}  imply
 \beg{align*} & \sup_{\|h_\nu\|_\infty\le C} \E^\nu|\Xi_r(t)-\Xi_r^{(n)}(t)| \\
 &=\ff 2 t  \sup_{\|h_\nu\|_\infty\le C} \sum_{k=n+1}^\infty \ff{\e^{-2\ll_k r} }{\ll_k^2} \int_0^t \nu(P_s \phi_k^2) (1-\e^{-\ll_k(t-s)})\d s \le C \vv_{n}, \end{align*}
 where    $\vv_{n}:= 2\sum_{k=n+1}^\infty \ff{2}{\ll_k^2\e^{2\ll_k r}}\to 0$  as $n\to\infty$.
   Combining this with \eqref{GN*1} we see that for any $a\in\R$ and $\vv>0$,
   \beg{align*} &\limsup_{t\to\infty}\sup_{\|h_\nu\|_\infty\le C} \big|\P^\nu(\Xi_r(t)<a)- \P(\Xi_r<a)\big|\\
   &\le \limsup_{t\to\infty}\sup_{\|h_\nu\|_\infty\le C} \Big\{\big|\P^\nu(\Xi_r^{(n)}(t)<a-\vv)- \P(\Xi_r^{(n)}<a-\vv)\big|\\
   &\qquad  \qquad + \big|\P^\nu(\Xi_r^{(n)}(t)<a-\vv)- \P^\nu(\Xi_r(t)<a)\big|\Big\}
    + \big|\P (\Xi_r^{(n)}<a-\vv)- \P(\Xi_r<a)\big|\\
   &\le  \limsup_{t\to\infty}\sup_{\|h_\nu\|_\infty\le C} \Big\{\P^\nu(|\Xi_r(t)-\Xi_r^{(n)}(t)|\ge \vv) +\P^\nu(a-\vv\le \Xi_r^{(n)}(t) <a)\Big\} \\
   &\qquad \qquad + \P(|\Xi_r-\Xi_r^{(n)}|\ge \vv) + \P(a-\vv\le \Xi_r^{(n)}<a) \\
   &\le \ff{(1+C)\vv_n}\vv + 2 \P(a-\vv\le \Xi_r^{(n)}<a),\ \ \vv>0, n\ge 1.\end{align*}
   Letting first $n\uparrow \infty$ then $\vv\downarrow 0$, we prove  \eqref{CLT1}.

(2)   Next, let  $M$ be compact with  $d\le 3$. We have $\sum_{k=1}^\infty \ff 2 {\ll_k^2}<\infty$, so that  the proof in (1) applies to $r=0$ or $r=r_t$ with $r_t\downarrow 0$ as $t\uparrow \infty$, where
$\Xi_0(t):= \sum_{k=1}^\infty \ll_k^{-1} |\psi_k(t)|^2, \Xi_0:= 2\sum_{k=1}^\infty \ll_k^{-2}\xi_k^2$. Then  \eqref{CLT2} holds.
\end{proof}

\beg{lem}\label{L02}  Assume {\bf (A1)}. For any $0<\vv<t,$ let
 $$\mu_{t,r}^\vv= \ff 1 {t-\vv} \int_{\vv}^t P_r(X_s,\cdot)\d s,\ \ r\ge 0,$$  where $P_r(X_s,\cdot):=\dd_{X_s}$ for $r=0.$ Let $D$ be the diameter of $M$. Then
\beg{align*} &|t\W_2(\mu_{t,r},\mu)^2- (t-\vv)\W_2(\mu_{t,r}^\vv,\mu)^2| \\
&\le 3c(r)\ss\vv +  \ss\vv (t-\vv) \W_2(\mu_{t,r}^\vv,\mu)^2,\ \ r\ge 0, t>\vv, \vv\in (0,1) \end{align*}
holds for
$c(r):=\min\big\{\|p_r\|_\infty^2 (\mu\times\mu)(\rr^2), D^2\big\},$ which is finite if either $r>0$ or $D<\infty$. \end{lem}
\beg{proof} It is easy to see that the measure
$$\pi(\d x,\d y):= \bigg(\ff 1 t\int_\vv^t P_r(X_s,\d x)\d s\bigg)\dd_x(\d y) + \bigg(\ff 1 {t(t-\vv)}\int_\vv^t P_r(X_s,\d x)\d s\bigg) \int_0^\vv P_r(X_s,\d y)\d s$$ is a coupling of $\mu_{t,r}^\vv$ and $\mu_t$. So,
\beg{align*} &t \W_2(\mu_{t,r}^\vv,\mu_{t,r})^2 \le t\int_{M\times M} \rr(x,y)^2 \pi(\d x,\d y) \\
&= \ff 1 {t-\vv} \int_\vv^t \d s_1 \int_0^\vv \d s_2\int_{M\times M}  \rr(x,y)^2 p_r(X_{s_1},x)p_r(X_{s_2}, y)\mu(\d x)\mu(\d y) \le c(r) \vv. \end{align*}
On the other hand,
$$\W_2(\mu_{t,r}^\vv,\mu)^2 \le \int_{M\times M} \rr(x,y)^2 \mu_{t,r}^\vv(\d x)\mu(\d y)\le c(r),\ \ r\ge 0.$$
Therefore,
\beg{align*} &|t\W_2(\mu_{t,r},\mu)^2- (t-\vv)\W_2(\mu_{t,r}^\vv,\mu)^2| \\
&\le \vv \W_2(\mu_{t,r},\mu)^2 +  (t-\vv) \big\{|\W_2(\mu_{t,r},\mu) -\W_2(\mu_{t,r}^\vv,\mu)|^2 + 2 |\W_2(\mu_{t,r},\mu) -\W_2(\mu_{t,r}^\vv,\mu)|  \W_2(\mu_{t,r}^\vv,\mu)\big\}\\
&\le \vv \W_2(\mu_{t,r},\mu)^2 + (1+\vv^{-\ff 1 2}) (t-\vv)|\W_2(\mu_{t,r},\mu) -\W_2(\mu_{t,r}^\vv,\mu)|^2+ \vv^{\ff 1 2} (t-\vv) \W_2(\mu_{t,r}^\vv,\mu)^2 \\
&\le 3c(r) \ss{ \vv} +  \ss \vv \,(t-\vv)\W_2(\mu_{t,r}^\vv, \mu)^2,\ \  t>\vv, \vv\in (0,1).\end{align*}
\end{proof}

\beg{proof}[Proof of \eqref{CLT}] (a) We first prove for $\nu$ with   $\|h_\nu\|_\infty\le C$.  Take $\si= t^{-\ff23}$. By \eqref{ANN} and \eqref{AN3},   we  find a constant $c_1>0$ such that for large enough $t\ge 1$,
 it holds  on the event $B_\si:= \{\|f_{t,r}-1\|_\infty\le\si^{\ff23} \}$ that
\beq\label{WMX} t\W_2(\mu_{t,r},\mu)^2 \ge \Xi_r(t) -c_1t\si^{\ff53}=  \Xi_r(t) -c_1t^{- \ff 1 9}. \end{equation}   Moreover, Lemma \ref{L2.3} with $k=1$ implies
\beq\label{*1*} \lim_{t\to\infty} \sup_{x\in M}\P^x (B_\si^c)\le  c(1,r)  \lim_{t\to\infty} \si^{-\ff43}   t^{-1} = c(1,r)  \lim_{t\to\infty}     t^{-\ff19}=0.\end{equation}
It follows from    \eqref{WMX} and \eqref{*1*} that
\beq\label{***} \lim_{t\to\infty} \sup_{x\in M} \P^x \big(t\W_2(\mu_{t,r}, \mu)^2\le (1-\vv) \Xi_r(t) -\vv\big) \le \lim_{t\to\infty} \sup_{x\in M} \P^x (B_\si^c) =0.\end{equation}
On the other hand, since $\scr M(r,1)\to 1$ as $r\to 1$, \eqref{*1*} implies that $\scr M(f_{t,r},1)\to 1$ in $\P^x$ uniformly in $x\in M$, so that \eqref{PP} implies
  $$\lim_{t\to\infty} \sup_{x\in M} \P^x\Big(t \W_2(\mu_{t,r},\mu)^2 \ge (1+\vv) \Xi_r(t) +\vv\Big) =0,\ \ \vv>0.$$ This together with \eqref{***} and \eqref{CC} yields
\beq\label{WM*0} \lim_{t\to\infty} \sup_{\|h_\nu\|_\infty\le C} \P^\nu(|t\W_2(\mu_{t,r},\mu)^2-\Xi_r(t)|\ge \vv)=0,\ \ \vv>0,C>0.\end{equation}  Combining this with \eqref{CLT1}  we prove
\beq\label{CLTT}  \lim_{t\to \infty}\sup_{\|h_\nu\|_\infty\le C } \big|\P^\nu\big(t\W_2(\mu_{t,r},\mu)^2 <a\big)-\nu_r\big((-\infty,a)\big)\big|=0,\ \ a\in \R.\end{equation}

(b) We now consider $\nu=\dd_x$.  By the Markov property, the law of $\mu_{t,r}^\vv $ under $\P^x$ coincides with that of $\mu_{t-\vv,r}$ under   $\P^\nu$ with $\nu(\d y):= p_\vv(x,y)\mu(\d y).$
Moreover, since $\sup_{x,y}p_\vv(x,y)=\|P_\vv\|_{1\to\infty}=:c(\vv)<\infty,$
   \eqref{CLTT}   implies
 $$\lim_{t\to\infty} \sup_{x\in M} \big|\P^x((t-\vv)\W_2(\mu_{t,r}^\vv,\mu)^2<a)- \nu_r((-\infty,a))\big|=0,\ \ a\in\R.$$ Combining this with Lemma \ref{L02}, we obtain
 \beg{align*} &\lim_{t\to \infty}\sup_{x\in M} \big|\P^x\big(t  \W_2(\mu_{t,r},\mu)^2 <a\big)-\nu_r\big((-\infty,a)\big)\big| \\
 &\le  \lim_{t\to \infty}\sup_{x\in M} \Big\{\big|\P^x\big(t  \W_2(\mu_{t,r}^\vv,\mu)^2 <a+3c(r)\ss\vv +\dd \big)-\nu_r\big((-\infty,a+3c(r)\ss\vv +\dd )\big)\big|\\
 &\qquad \qquad+ \big|\P^x\big(t  \W_2(\mu_{t,r}^\vv,\mu)^2 <a-3c(r)\ss\vv +\dd  \big)-\nu_r\big((-\infty,a-3c(r)\ss\vv +\dd )\big)\big|\\
  &\qquad\qquad + \P^x\big(\ss\vv (t-\vv) \W_2(\mu_{t,r}^\vv,\mu)^2\ge \dd\big)\Big\}+  \nu_r([a-3c(r)\ss\vv +\dd, a+3c(r)\ss\vv +\dd])\\
&= \nu([\dd\vv^{-\ff 1 2},\infty))  +\nu_r([a-3c(r)\ss\vv +\dd, a+3c(r)\ss\vv +\dd]),\ \ \vv,\dd >0.\end{align*} Taking $\dd= \vv^{\ff 1 4}$ and letting $\vv\to 0$ we finish the proof.
\end{proof}

\beg{proof}[Proof of Corollary \ref{C1.2}] Obviously, when $d\le 3$, \eqref{E2'} and \eqref{*2}   imply assertion (1). Next, for $d=4$, \eqref{*2} implies
\beq\label{N*1} c_1' \sum_{i=1}^\infty i^{-1} \e^{-c_2' r\ss i} \le \sum_{i=1}^\infty \ff 2 {\ll_i^2\e^{2r\ll_i}}\le c_1 \sum_{i=1}^\infty i^{-1} \e^{-c_2 r\ss i},\ \ r>0\end{equation}
for some constants $c_1,c_2,c_1',c_2'>0.$ Moreover, there exist constants $c_3,c_4>0$ such that
 \beq\label{N**1}\beg{split} &\sum_{i=1}^\infty i^{-1} \e^{-c_2 r\ss i}\le c_3\int_1^\infty s^{-1}\e^{-\ff{c_2r}2\ss s}\d s \\
&=c_3 \int_{r^2}^\infty t^{-1} \e^{-\ff{c_2}2 \ss t}\d t \le c_4 \log r^{-1},\ \ r\in (0,1/2),\end{split}\end{equation}  while
for some constants $c_3',c_4'>0$,
\beq\label{OPL} \beg{split}  &\sum_{i=1}^\infty i^{-1} \e^{-c_2' r\ss i}\ge c_3'\int_1^\infty s^{-1}\e^{-c_2'r\ss s}\d s \\
&=c_3' \int_{r^2}^\infty t^{-1} \e^{ -c_2' \ss t}\d t \ge c_4' \log r^{-1},\ \ r\in (0,1/2).\end{split}\end{equation} Combining this with \eqref{N*1}, \eqref{N**1} and \eqref{E2'}, we prove the second assertion.

Finally, when $d\ge 5$, \eqref{*2} implies that for some constants $c_i,c_i', i=1,2,3$ such that
\beq\label{N*2} \beg{split} &\sum_{i=1}^\infty \ff 2 {\ll_i^2\e^{2r\ll_i}}\le c_1 \sum_{i=1}^\infty i^{-\ff 4 d} \e^{-c_2 r i^{\ff 2 d}} \le c_1 \int_0^\infty s^{-\ff 4 d} \e^{-c_2 r s^{\ff 2 d}}\d s \\
& = c_1 \int_0^\infty r^{\ff{4-d}2} t^{-\ff 4 d}\e^{-c_2t^{\ff 2 d}}\d t \le c_3 r^{\ff{4-d}2},\ \ r>0,\end{split}\end{equation} and
\beq\label{N*3} \beg{split} &\sum_{i=1}^\infty \ff 2 {\ll_i^2\e^{2r\ll_i}}\ge c_1' \sum_{i=1}^\infty i^{-\ff 4 d} \e^{-c_2' r i^{\ff 2 d}} \ge c_1 \int_1^\infty s^{-\ff 4 d} \e^{-c_2' r s^{\ff 2 d}}\d s \\
& = c_1' \int_{r^{\ff d 2}}^\infty r^{\ff{4-d}2} t^{-\ff 4 d}\e^{-c_2't^{\ff 2 d}}\d t \ge c_3' r^{\ff{4-d}2},\ \ r\in(0,1),\end{split}\end{equation}
Combining these with \eqref{E2'}, we prove (3). \end{proof}

\section{Proof of Theorem \ref{T1.1}}
In this section we assume that $M$ is compact. We first present some lemmas which will be used in the proof.

\subsection{Some lemmas}

When $M$ is compact, there exist constants $\kk,\ll>0$ such that
  \beq\label{IIT} \|P_t\|_{p\to q} \le \kk  (1\land t)^{-\ff d 2(p^{-1}-q^{-1})},\ \ t>0, q\ge p\ge 1.\end{equation}
 In particular, {\bf (A1)} holds with $\|P_t\|_{1\to\infty}\le \kk (1\land t)^{-\ff d 2}$ for some constant $\kk>0$ and all $t>0$, so that \eqref{E1'} follows from Theorem \ref{T1.2}.

To estimate $\E[W_2(\mu_t,\mu)^2]$ from \eqref{E1'},   we use the triangle inequality to derive
\beq\label{TRA} \E[\W_2(\mu_t,\mu)^2]\le (1+\vv) \E[\W_2(\mu_{t,r},\mu)^2] +(1+\vv^{-1}) \E[\W_2(\mu_t,\mu_{t,r})^2],\ \ \vv>0.\end{equation}
We will show that $\E[\W_2(\mu_t,\mu_{t,r})^2]\le c r$ holds for some constant $c>0$ and all $r>0$, which is   known when $\pp M$ is either empty or convex, but is new when $\pp M$ is non-convex, see \eqref{LK1} below.  If we could take $ r_t>0$    such that
$$\lim_{t\to\infty} tr_t=0,\ \ \limsup_{t\to\infty} \big\{t\E\W_2(\mu_{t,r_t},\mu)^2\big\}\le \sum_{i=1}^\infty\ff 2 {\ll_i^2},$$ we would  deduce  the desired estimate \eqref{E1} from \eqref{TRA}.
Let us start with the following estimate of $\E[\W_2(\mu_t,\mu_{t,r})^2]$.

\beg{lem}\label{LN4}  Assume that $M$ is compact and let $\mu_{t,r,\vv}= (1 - \vv)\mu_{t, r}+\vv\mu, \vv\in [0,1].$  Then there exists a constant $c>0$  such that
\beq\label{LK1}   \E^\nu [\W_2(\mu_t,\mu_{t, r})^2]\le c \|h_\nu\|_\infty r,\ \  \nu= h_\nu\mu, r\ge 0,\end{equation} and for any initial value $X_0$ of the diffusion process,
\beq\label{LK2} \W_2(\mu_{t,r,\vv},\mu_{t,r})^2\le c\vv, \ \  t,r>0, \vv\in [0,1].\end{equation} \end{lem}

\beg{proof} Since $M$ is compact, there exists $\tt\rr\in C_b^\infty(M\times M)$ and constants $\aa_2\ge \aa_1>0$ such that
 $$\aa_1 \tt\rr \le \rr \le \aa_2 \tt\rr.$$
 By It\^o's formula,     there exist  constants $c_1,c_2>0$ such that
   \beq\label{POI}  \beg{split} \d \tt\rr(X_0,X_r)^2&=\big\{ L \tt\rr (X_0,\cdot)^2(X_r)\big\}\d r + \d M_r + \big\{N\tt\rr(X_0,\cdot)^2(X_r)\big\}\d l_r\\
 &\le c_1 \d r + \d M_r + c_2 \d l_r,\end{split}\end{equation}
 where $M_r$ is a martingale,  when $\pp M$ exists $N$ is the inward unit normal vector field of $\pp M$ and $l_r$ is the local time of $X_r$ on $\pp M$, and $D$ is the diameter of $M$. If $\pp M =\emptyset$, then $l_r=0$  so that
\beq\label{ACC}  \E^\nu \big[\rr(X_0,X_r)^2\big]\le \aa_2^2 \E^\nu \big[\tt\rr(X_0,X_r)^2\big] \le c_1\aa_2^2 r\le c_1\aa_2^2\|h_\nu\|_\infty r,\ \ \ r\ge 0.\end{equation}
 When $\pp M\ne \emptyset$, \eqref{POI} implies
 \beq\label{KK} \E^\nu\big[ \rr(X_0,X_r)^2\big]\le c_1\aa_2^2 r + c_2\aa_2^2 \E^\nu l_r,\ \ r>0.\end{equation}
 Let $\tau=\inf\{t\ge 0: X_t\in \pp M\}$. We have  $l_r=0$ for $r\le \tau$, so that  by the Markov property
\beq\label{KK2} \E^\nu l_r = \E^\nu [1_{\{\tau<r\}}\E^{X_\tau} l_{r-\tau} ] \le \P^\nu(\tau<r) \sup_{x\in \pp M} \E^x l_r.\end{equation}
 By  \cite[Proposition 4.1]{W10} and \cite[Lemma 2.3]{ATW09}, there exist constants $c_2,c_3,c_4>0$ such that
 \beg{align*} &\E^xl_r\le c_2 \ss r, \ \ x\in\pp M, \\
 &\P^\nu(\tau<r)\le \int_M \e^{-c_2\rr_\pp(x)^2/r}\nu(\d x) \le \|h_\nu\|_\infty \int_M \e^{-c_3\rr_\pp(x)^2/r}\mu(\d x)\le c_4\|h_\nu\|_\infty \ss r.\end{align*}
Combining these with  \eqref{KK2} we derive  $\E^\nu l_r\le c_2c_4 \|h_\nu\|_\infty r$ for $r\ge 0.$ Therefore,  by \eqref{KK} for $\pp M\ne \emptyset$ and    \eqref{ACC} for $\pp M=\emptyset$,  we  find a constant $c>0$ such that in any case
 \beq\label{ACC1}   \E^\nu \big[\rr(X_0,X_r)^2\big]\le c \|h_\nu\|_\infty  r,\ \ \ r\ge 0. \end{equation}
It is easy to see that for any $t>0$,
 $$\pi_t(\d x,\d y):= \bigg(\ff 1 t \int_0^t\big\{p_r(x,y) \dd_{X_s} \big\}(\d x)  \d s\bigg)\mu(\d y) \in \scr C(\mu_t,\mu_{t,r}).$$
Then
  \beq\label{ACC2} \beg{split} & \W_2(\mu_{t,r}, \mu_t)^2 \le \iint_{M \times M} \rr(x,y)^2 \pi_t(\d x,\d y)\\
   &=\ff 1 t \int_0^t \d s \int_M p_r(X_s,y) \rr(X_s,y)^2 \mu(\d y),\ \ r,t>0. \end{split}\end{equation}
  Letting $\nu_s= (P_sh_\nu)\mu$, which is the distribution of $X_s$ provided the law of $X_0$ is $\nu$, by the Markov property and \eqref{ACC1}, we obtain
  $$\E^\nu\int_M p_r(X_s,y)\rr(X_s,y)^2\mu(\d y)= \E^{\nu_s} \big[\rr(X_0,X_r)^2\big] \le c\|P_sh_\nu\|_\infty r \le c \|h_\nu\|_\infty r,\ \ s,r>0.$$
  Substituting this into \eqref{ACC2}, we prove    \eqref{LK1}.

  On the other hand, since $\mu_{t,r,\vv}=(1-\vv)\mu_{t,r}+\vv\mu$,  we have $$\pi (\d x,\d y):= (1-\vv) \mu_{t,r}(\d x) \dd_x(\d y)+ \vv \mu(\d x) \mu_{t,r}(\d y)\in \C(\mu_{t,r,\vv},\mu_{t,r}),$$ so that
  $$\W_2(\mu_{t,r,\vv}, \mu_{t,r})^2\le \int_{M\times M}\rr(x,y)^2 \pi (\d x,\d y)\le \vv D^2.$$
  Therefore, \eqref{LK2} holds.
\end{proof}

The following lemmas is used for the estimates of $\E[\W_2(\mu_{t,r},\mu)^2]$. Let us start with the case $d \le 3$.

\beg{lem}\label{L2.3N}  Assume that $M$ is compact with $d \le 3$.   There exists a constant $c > 0$ such that
for any probability measure $\nu=h_\nu\mu$,
\beq\label{L2.31} \sup_{y\in M} \E^\nu\big[|f_{t,r}(y)-1|^2\big] \le    \ff{ c \|h_\nu\|_\infty}{t\ss r}   , \, \ t\ge 1, r \in (0, 1], \end{equation}
\end{lem}

\beg{proof} We use the notation in the proof of Lemma \ref{L2.3}.   Noting that $f = p_r(\cdot,y)-1$ and $M$ is compact, by \eqref{SPP} and \eqref{IIT} there exists a constant $c > 0$ such that 
 \beg{align*}&\E^\nu |g(r_1,r_2)| = \nu(P_{r_1} |fP_{r_2-r_1}f|)\le \|h_\nu\|_{\infty} \mu(|fP_{r_2-r_1}f|) \\
  &\le  2 \|h_\nu\|_{\infty} \| P_{\ff{r + r_2 - r_1} 2 }P_{\ff {r_2 - r_1} 2}(p_{r/2}(\cdot,y)-1)\|_\infty\\
  &\le c \|h_\nu\|_{\infty} (1 \land (r+r_2-r_1))^{-\ff{d}2} \e^{- \ff{\ll_1}2 (r_2-r_1)}.\end{align*}  So,
\beq\label{OOO} \beg{split} & \int_{\DD_1(t)} \E^\nu|g(r_1,r_2)| \d r_1\d r_2
\le  c \|h_\nu\|_\infty
\int_0^t   \d r_1 \int_{r_1}^t (1 \land (r+r_2-r_1))^{-\ff{d}2}\e^{- \ff{\ll_1}2 (r_2-r_1)}\d r_2\\
& \le  c \|h_\nu\|_\infty
\int_0^t   \d r_1 \int_{r_1}^t \big(1 + (r+r_2-r_1)^{-\ff{d}2} \big) \e^{- \ff{\ll_1}2 (r_2-r_1)}\d r_2.  \end{split} \end{equation}
Noting that
$$\int_{r_1}^t (r+r_2-r_1)^{-\ff{d}{2} }\e^{- \ff{\ll_1}2 (r_2-r_1)}\d r_2\le \int_0^\infty (r+s)^{-\ff{d}{2} }\e^{-\ff{\ll_1}2 s}\d s\le c' \beg{cases} 1,&\text{when}\, \ d = 1   ; \\ \log(1+ r^{-1} ),&\text{when}\, \ d = 2 ;\\ r^{1 - \ff d 2 },&\text{when}\, \ d \ge 3,   \end{cases} $$ 
holds  for some constant $c'>0$, combining \eqref{**1} with \eqref{OOO} for $k = 1$ and noticing that $d \le 3$, we prove \eqref{L2.31}.
\end{proof}

\beg{lem}\label{LN1}  Assume that $M$ is compact with $d \le 3$. For any $\aa\in (1,2)$ and $r_t:=t^{-\aa}$,
\beq\label{*HU1} \lim_{t\to\infty } \sup_{\|h_\nu\|_\infty\le C, y\in M} \E^\nu \big[\big|\M((1-r_t)f_{t,r_t}(y)+r_t,1)^{-1}-1\big|^q\big]=0,\ \ C,q>0.\end{equation}
\end{lem}
\beg{proof} By \cite[Lemma 3.12]{AMB},
\beq\label{*00}  \ff{\theta (ab)^{\ff \theta 2}|a-b|}{|a^\theta-b^\theta|}\le \M(a,b)\le \ff{\theta(a^\theta+b^\theta)(a-b)}{2(a^\theta-b^\theta)},\ \ a,b,\theta >0.\end{equation}
Combining this with the simple inequality $|a^\theta-1|\le |a-1|$ for $a\ge 0$ and $\theta\in [0,1]$, we obtain
 \beg{align*} &\M((1-r)f_{t,r}(y)+r,1)\ge \ff{\theta
\{(1-r)f_{t,r}(y)+r\}^{\ff\theta 2}|(1-r)f_{t,r}(y)+r-1|}{|\{(1-r)f_{t,r}(y)+r\}^\theta-1|}\\
&\ge \theta \{(1-r)f_{t,r}(y)+r\}^{\ff \theta 2}\ge \theta r^{\ff\theta 2},\ \ t\ge 1, \theta\in (0,1), r>0.\end{align*}
This implies
\beq\label{*01} \big|\M((1-r)f_{t,r}(y)+r,1)^{-1}-1\big|\le 1 + \theta^{-1} r^{-\ff{\theta} 2},\ \ t\ge 1, \theta\in (0,1), r>0.\end{equation}

On the other hand, let $\eta\in (0,1)$. On the event
$$A_{\eta,y}:= \{|f_{t,r}(y)-1|\le\eta\}$$ we have
$|(1-r)f_{t,r}(y)+r -1|\le\eta,$ so that \eqref{*00} for $\theta=1$ implies
$$\ss{1-\eta}\le\M((1-r)f_{t,r}(y)+r,1) \le 1+\ff \eta 2\ \text{on}\ A_{\eta,y}.$$
Thus,
$$1_{A_{\eta,y}} \big|\M((1-r)f_{t,r}(y)+r,1)^{-1}-1\big|^q\le \Big|\ff 1{\ss{1-\eta}}- \ff 2 {2+\eta}\Big|^q=:\dd_\eta.$$
  Combining this with \eqref{L2.31} and using \eqref{*01}, we obtain that for $t \ge 1$ and $r \in (0, 1]$,
\beq\label{*HU} \beg{split} &\sup_{y\in M}\E^\nu \big[\big|\M((1-r)f_{t,r}(y)+r,1)^{-1}-1\big|^q\big]\le \dd_\eta + (1+\theta^{-1}r^{-\ff{\theta} 2})^q\sup_{y\in M} \P^\nu(A_{\eta,y}^c)
\\
&\le  \dd_\eta + C(\theta,\eta)   \|h_\nu\|_\infty t^{-1}
r^{- \ff 1 2 -\ff {\theta q}{2}}.  \end{split} \end{equation}
Then for any $\aa\in (1,2)$ and $q>0$, we may take a small enough $\theta$ such that
$\aa(\ff 1 2 +\ff{\theta q}2) < 1.$ Then \eqref{*HU1} follows from \eqref{*HU} with $r = r_t = t^{-\aa}$ and $\eta\downarrow 0$.
\end{proof}

\beg{lem}\label{LN2}  Assume that $M$ is compact. For any $p\in [1,2]$, there exists a constant
$c>0$ such that $\psi_i(t):=\ff 1 {\ss t}\int_0^t \phi_i(X_s)\d s$ satisfies
$$\E^\nu\big[ |\psi_i(t)|^{2p}\big]\le  c \|h_\nu\|_\infty    \ll_i^{p-2+(p-1)(\ff d 2 -2)},\ \ t\ge 1, i\in\mathbb N,\nu=h_\nu\mu.$$\end{lem}

\beg{proof} Let $f=\phi_i$. Then $g(r_1,r_2)$ in \eqref{**1} satisfies
\beq\label{Q1} g(r_1,r_2)=(\phi_iP_{r_2-r_1}\phi_i)(X_{r_1})=\e^{-(r_2-r_1)\ll_i}\phi_i(X_{r_1})^2.\end{equation}
Since $\mu(h_\nu P_{r_1}\phi_i^2)\le \|h_\nu\|_\infty\mu(\phi_i^2)=\|h_\nu\|_\infty<\infty,$ this and  \eqref{**1} with $k=1$ imply
 \beq\label{Q2} \beg{split} &t\E^\nu \big[|\psi_i(t)|^2\big] \le c_1  \int_0^t \d r_1 \int_{r_1}^t \E^\nu [g(r_1,r_2)]\,\d r_2 \\
&= c_1 \int_0^t \d r_1 \int_{r_1}^t \e^{-(r_2-r_1)\ll_i}\mu(h_\nu P_{r_1}\phi_i^2)\,\d r_2\le c_1 \|h_\nu\|_\infty \ff t{\ll_i},\ \ t\ge 1, i\in \mathbb N \end{split}\end{equation} for some constant $c_1>0$.
On the other hand, taking $k=2$ in \eqref{**1} and using \eqref{Q1}, we find a constant $c_2>0$ such that
\beg{align*}&t^2 \E^\nu \big[|\psi_i(t)|^4\big]\le c_2 \bigg(\int_0^t\d r_1\int_{r_1}^t \big(\E^\nu |g(r_1,r_2)|^2\big)^{\ff 1 2} \d r_2\bigg)^2 \\
&=c_2 \bigg(\int_0^t\d r_1\int_{r_1}^t \e^{-(r_2-r_1)\ll_i} \ss{\mu(h_\nu P_{r_1}\phi_i^4)} \,   \d r_2\bigg)^2,\ \  t\ge 1,i\in\mathbb N.\end{align*}
By \eqref{IIT} and $P_t \phi_i=\e^{-\ll_i t}\phi_i$, we obtain
\beq\label{HTQ} \|\phi_i\|_\infty=\inf_{t>0} \big\{\e^{\ll_i t}\|P_t\phi_i\|_\infty\big\}\le \inf_{t>0} \big\{\e^{\ll_i t}\|P_t\|_{2\to\infty} \big\}  \le c_3\ll_i^{\ff d 4},\ \ i\ge 1\end{equation} for some constant $c_3>0$.
Since $h_\nu$ is bounded,  \eqref{HTQ} and $\mu(\phi_i^2)=1$ imply
$$\ss{\mu(h_\nu P_{r_1} \phi_i^4)}\le \ss{\|h_\nu\|_\infty \mu(\phi_i^4)} \le \ss{\|h_\nu\|_\infty \|\phi_i\|_\infty^2\mu(\phi_i^2)}\le  c_3 \ss{ \|h_\nu\|_\infty}\,  \ll_i^{\ff d 4},\ \ i\ge 1.$$
Therefore, there exists a constant $c_4>0$ such that
$$ t^2\E^\nu \big[|\psi_i(t)|^4\big]  \le c_4 \|h_\nu\|_\infty t^2\ll_i^{\ff d 2-2 },\ \  t\ge 1,i\in\mathbb N.$$
Combining this with \eqref{Q2} and H\"older's inequality,  we find a constant $c>0$ such that for any $p\in [1,2]$,
\beg{align*} &\E^\nu \big[|\psi_i(t)|^{2p}\big] =\E^\nu \big[|\psi_i(t)|^{4-2p}|\psi_i(t)|^{4(p-1)}\big]\\
&\le \big(\E^\nu |\psi_i(t)|^2\big)^{2-p} \big(\E^\nu |\psi_i(t)|^4\big)^{p-1}
 \le  c\|h_\nu\|_\infty  \ll_i^{p-2 +(p-1) (\ff d 2 -2)},\ \ t\ge 1, i\in \mathbb N.\end{align*}

\end{proof}

\beg{lem}\label{LN3}  Assume that $M$ is compact with $d \le 3$.
There exists a constant $p>1$ such that for any $C>1$,
$$\limsup_{t\to\infty}   \sup_{r>0,\|h_\nu\|_\infty\le C}   \bigg\{t^p\,\E^\nu\int_M |\nn L^{-1} (f_{t,r}-1)|^{2p}\d\mu\bigg\} <\infty.$$
\end{lem}

\beg{proof} Let $p>1$. By \cite[(1.10)]{W05}, the gradient estimate
\beq\label{B0} |\nn P_t f|\le \ff {c(p)}{\ss t} (P_t|f|^p)^{\ff 1 p},\ \ t>0, f\in \B_b(M)\end{equation}
holds for some constant $c(p)>0$. Combining this with
  \eqref{LLL} and \eqref{FT},    we obtain
 \beq\label{90} \beg{split} &\E^\nu \int_M |\nn L^{-1} (f_{t,r}-1)|^{2p}\d\mu\le   \E^\nu \int_M \bigg(\int_0^\infty |\nn P_s (f_{t,r}-1)|\d s\bigg)^{2p}\d\mu\\
&\le c_1(p) \E^\nu \int_M \bigg(\int_0^\infty \ff 1 {\ss s} \big\{P_{\ff s 2}| P_{\ff s 2}(f_{t,r}-1)|^p\big\}^{\ff 1 p} \d s\bigg)^{2p}\d\mu\\
&\le c_1(p)\ \bigg(\int_0^\infty s^{-\ff{2p}{2(2p-1)}} \e^{-\ff{2p\theta s}{2p-1}} \d s\bigg)^{2p-1}\\
&\qquad\qquad\times
\E^\nu   \int_0^\infty \e^{\theta s}   \mu\big( \{P_{\ff s 2}| P_{\ff s 2}f_{t,r}-1|^p \}^{2}\big) \d s,\ \ t\ge 1, r>0 \end{split}\end{equation} for some constant $c_1(p)>0$.
Let $\theta\in (0,\ff{\ll_1}2)$ and $p\in (1,2)$. We have
\beq\label{LOU} \int_0^\infty s^{-\ff{2p}{2(2p-1)}} \e^{-\ff{2p\theta s}{2p-1}} \d s<\infty.\end{equation}
Combining \eqref{NNK}, \eqref{IIT}, \eqref{90}, \eqref{LOU} and H\"older's inequality,  we arrive at
  \beg{align*} &t^p \E^\nu \int_M |\nn L^{-1} (f_{t,r}-1)|^{2p}\d\mu
 \le c_2(p) t^p \E^\nu \int_0^\infty \e^{\theta s} \|P_{\ff s 2}\|_{\ff 2 p\to 2}^2\big\{\mu((P_{\ff s 2 } (f_{t,r}-1))^2)\big\}^p\d s  \\
&\le c_3(p)  \E^\nu\int_0^\infty \e^{\theta s} (1\land s)^{-\ff{d(p-1)}2} \Big( \sum_{i=1}^\infty \e^{-(2r+s)\ll_i} |\psi_i(t)|^{2} \Big)^p\d s\\
&\le c_3(p)  \Big(\sum_{i=1}^\infty  i^{-\ff{\vv}{p-1}}\Big)^{p-1}   \int_0^\infty   (1\land s)^{-\ff{d(p-1)}2}  \sum_{i=1}^\infty i^\vv \e^{-p(2r+s)\ll_i+\theta s}
\E^\nu\big[|\psi_i(t)|^{2p}\big]  \d s,\ \ t\ge 1,i\in\mathbb N \end{align*}  for some constants $c_2(p), c_3(p)>0$. Since
$-ps\ll_i+\theta s\le -\ff s 2 \ll_i,$ and noting that  for any $c>0$ and $\dd\in (0,1)$ there exists a constant $c'>0$ such that
$$\int_0^\infty (1\land s)^{-\dd}\e^{-c\ll_i s}\d s\le c' \ll_i^{\dd-1},\ \  i\ge 1,$$
this implies
$$t^p\E^\nu \int_M |\nn L^{-1} (f_{t,r}-1)|^{2p}\d\mu \le c_4(p)  \Big(\sum_{i=1}^\infty  i^{-\ff{\vv}{p-1}}\Big)^{p-1}    \sum_{i=1}^\infty i^\vv  \ll_i^{\ff{d(p-1)}2-1}\e^{-2r\ll_i}  \E^\nu\big[|\psi_i(t)|^{2p}\big]$$ for some constant $c_4(p)>0$. Therefore,
 for any $\vv>0$ and $p>1$ such that  $\ff{\vv}{p-1}>1$, there exists  a constant $c(p,\vv)>0$ such that   this, \eqref{*2}  and Lemma \ref{LN2} yield
\beq\label{OOP} \beg{split} &\E^\nu \int_M |\nn L^{-1} (f_{t,r}-1)|^{2p}\d\mu\le c(p,\vv)t^{-p} \|h_\nu\|_\infty \sum_{i=1}^\infty i^{\dd_{p,\vv}}  \e^{-2ri^{2/d}},\ \ t\ge 1,r>0,\\
&\dd_{p,\vv}:= \vv+ \ff 2 d\big\{ (p-1)(d-2)+ p-3\big\}.\end{split} \end{equation}

When $d \le 3$,  by taking for instance   $\vv =\ff 1 {12}$, and   $p>1$ close enough to $1$ such that
  \beq\label{INF}\ff{d(p-1)}2 \in(0, 1),\ \  \ff{\vv}{p-1}>1, \ \   (p-1)(d-2)    -1 + p-2   \le -\ff 7 4,\end{equation}  and noting   $d\le 3$ and  \eqref{*2} imply   $\ll_i\ge c'' i^{\ff 2 3}$ for some constant $c''>0$,   from  \eqref{OOP}  we find a constant $c>0$ such that
\beg{align*} & \E^\nu \int_M |\nn L^{-1} (f_{t,r}-1)|^{2p}\d\mu \le c t^{-p} \|h_\nu\|_\infty \sum_{i=1}^\infty i^{\vv-\ff 2 3\cdot\ff 7 4} \\
&= c t^{-p} \|h_\nu\|_\infty \sum_{i=1}^\infty i^{\vv -\ff {7}{6} } <\infty,\ \ t\ge 1, r>0, \nu=h_\nu\mu.  \end{align*}
We finish the proof.
\end{proof}

For $d\ge 4$ we will use the following lemma  taken from \cite{L17}.  Although this lemma  is less sharper than Lemma \ref{L2.1}, it is easier to apply since the term $\M(1, f)$ is dropped.

\beg{lem} \label{L2.1'} Let $f$ be probability density functions with respect to $\mu$.   Then
 $$\W_2(f\mu, \mu)^2 \le 4 \int_M |\nn L^{-1}(f - 1)|^2 \d\mu.$$
 \end{lem}

Finally,  we have the following result on the large deviation of the empirical measures.
Let $\scr P$ be the set of all probability measures on $M$.

\beg{lem}[\cite{Wu2}]\label{LDP} Let $M$ be compact. Then for any open set $G\subset \scr P$ and closed set $F\subset \scr P$ under the weak topology,
 \beg{align*} &-\inf_{\nu\in G} I_\mu(\nu)\le \liminf_{t\to\infty} \ff 1 t \log  \inf_{x\in M} \P^x(\mu_t\in G), \\
& -\inf_{\nu\in F} I_\mu(\nu)\ge  \limsup_{t\to\infty} \ff 1 t \log  \sup_{x\in M} \P^x(\mu_t\in F). \end{align*}  \end{lem}

\beg{proof} Since the $\tau$-topology induced by bounded measurable functions on $M$ is stronger than the weak topology, $G$ and $F$ are open and closed respectively under the $\tau$-topology.
 By the ultracontractivity and irreducibility of $P_t$, \cite[Theorem 5.1(b) and Corollary B.11]{Wu2} imply
  \beg{align*} &-\inf_{\nu\in G} I_\mu(\nu)\le \liminf_{t\to\infty} \ff 1 t \log  {\rm ess}_\mu\inf_{x\in M} \P^x(\mu_t\in G), \\
& -\inf_{\nu\in F} I_\mu(\nu)\ge  \limsup_{t\to\infty} \ff 1 t \log  {\rm ess}_\mu\sup_{x\in M} \P^x(\mu_t\in F). \end{align*}  Consequently, letting $\scr P_c$ be the set of all probability measures on $M$
which are absolutely continuous with respect to $\mu$, we have
\beq\label{RPP} \beg{split}  &-\inf_{\nu\in G} I_\mu(\nu)\le \liminf_{t\to\infty} \ff 1 t \log   \inf_{\nu\in \scr P_c} \P^\nu(\mu_t\in G), \\
& -\inf_{\nu\in F} I_\mu(\nu)\ge  \limsup_{t\to\infty} \ff 1 t \log   \sup_{\nu\in \scr P_c} \P^\nu(\mu_t\in F). \end{split} \end{equation}
To replace $\scr P_c$ by $\scr P$, consider $\tt \mu_t^\vv:=\ff 1 {t}\int_\vv^{t +\vv}\dd_{X_s}\d s$ for $\vv>0$. By the Markov property, the law of $\tt\mu_t^\vv$ under $\P^x$ coincides with that of $\mu_t$ under $\P^\nu$,
where $\nu:=p_\vv(x,\cdot)\mu\in \scr P_c$.  So, \eqref{RPP} implies
\beq\label{RPP1} \beg{split}& -\inf_{\nu\in G} I_\mu(\nu)\le \liminf_{t\to\infty} \ff 1 t \log   \inf_{x\in M} \P^x(\tt\mu_t^\vv\in G)  \\
&-\inf_{\nu\in F} I_\mu(\nu)\ge  \limsup_{t\to\infty} \ff 1 t \log   \sup_{x\in M} \P^x (\tt\mu_t^\vv\in F),\ \ \vv>0. \end{split} \end{equation}

(a)  Let $D$ be the diameter of $M$. By taking the Wasserstein coupling
$$\pi(\d x,\d y):= (\mu_t\land \tt\mu_t^\vv)(\d x)\dd_x(\d y)+ \ff{(\mu_t-\tt\mu_t^\vv)^+(\d x)(\mu_t-\tt\mu_t^\vv)^-(\d y)}{(\mu_t-\tt\mu_t^\vv)^+(M)}\in \scr C(\mu_t,\tt\mu_t^\vv),$$
we obtain
\beq\label{RPP2} \W_2(\mu_t,\tt\mu_t^\vv)^2 \le \int_{M\times M}\rr^2\d\pi\le D^2 (\mu_t-\mu_t^\vv)^+(M)= \ff {D^2\vv}t.\end{equation}
So, when $t>D^2$, we have $\{\tt\mu_t^\vv\in G_\vv\}\subset \{\mu_t\in D\},$ where
$$G_\vv:=\Big\{\nu\in \scr P: \W_2(\nu, G^c)^2:=\inf_{\nu'\in G^c} \W_2(\nu, \nu')^2>\vv\Big\}$$
is an open subset of $\scr P$ under the weak topology, since for compact $M$,  $\W_2$  is continuous under the weak topology.
Combining this with \eqref{RPP1} for $G_\vv$ replacing $G$, we arrive at
$$-\inf_{\nu\in G_\vv} I_\mu(\nu)\le \liminf_{t\to\infty} \ff 1 t \log   \inf_{x\in M} \P^x(\tt\mu_t^\vv\in G),\ \ \vv>0.$$
Noting that $G_\vv\uparrow G$ as $\vv\downarrow 0$, we have $\inf_{\nu\in G_\vv} I_\mu(\nu)\downarrow  \inf_{\nu\in G} I_\mu(\nu)$ as $\vv \downarrow 0$. So, letting $\vv\downarrow 0$ we prove the desired inequality for open $G$.

(b) Similarly,   let $F_\vv:=\{\nu\in\scr P: \W_2(\nu, F)^2:=\inf_{\nu'\in F} \W_2(\nu,\nu')^2 \le \vv\}$, which is closed. When $t>D^2,$ \eqref{RPP1}  and \eqref{RPP2} imply
$$-\inf_{\nu\in F_\vv} I_\mu(\nu)\ge \limsup_{t\to\infty} \ff 1 t \log   \sup_{x\in M} \P^x(\mu_t\in F),\ \ \vv>0.$$ So, it suffices to show that
\beq\label{ALL} c:=\lim_{\vv\downarrow 0} \inf_{ F_\vv} I_\mu = \inf_{ F} I_\mu.\end{equation} Since $F_\vv\downarrow F$ as $\vv\downarrow 0$, we have   $c\le  \inf_{ F} I_\mu.$ On the other hand,
if $c<\infty$, then we may choose $\vv_n\downarrow 0$ and  $\nu_n\in F_{\vv_n}$ such that
  \beq\label{POU} I_\mu(\nu_n)\le \inf_{F_{\vv_n}} I_\mu+\ff 1 n \le \inf_{F_1}I_\mu+1<\infty,\ \ n\ge 1.\end{equation} So, $\nu_n= f_n\mu$ with $\sup_{n\ge 1}\mu(|\nn f_n^{\ff 1 2}|^2)<\infty$. By the Sobolev embedding theorem, $\{f_n^{\ff 1 2}\}_{n\ge 1}$ is relatively compact in $L^2(\mu)$, so that up to a subsequence
$f_n^{\ff 1 2}\to f^{\ff 1 2}$ in $L^2(\mu)$ for some probability density $f$ with respect to $\mu$. This and \eqref{POU} yield
$f^{\ff 1 2}\in W^{2,1}(\mu)$ and
$$I_\mu(f\mu):=\mu(|\nn f^{\ff 1 2}|^2) \le \liminf_{n\to\infty}I_\mu(\nu_n)\le \liminf_{n\to\infty} \inf_{F_{\vv_n}}I_\mu=c.$$
Since $F$ is closed,   $\nu_n\in F_{\vv_n}\downarrow F$ and $\nu_n\to f\mu$ weakly as $n\uparrow \infty$,  we conclude that $f\mu\in F$. Therefore, $\inf_FI_\mu\le c$ as desired.

\end{proof}

\subsection{Proof of Theorem \ref{T1.1} }

\beg{proof}[Proof of Theorem \ref{T1.1}(1)]
Obviously, \eqref{ULT} can be reformulated as  the following two estimates:
\beq\label{E1}  \limsup_{t\to \infty} \sup_{x\in M} \Big\{t \E^x [\W_2(\mu_{t},\mu)^2]\Big\}\le   \sum_{i=1}^\infty\ff 2 {\ll_i^2},  \end{equation}
\beq\label{E2} \liminf_{t\to\infty} \Big\{t\inf_{x\in M} \E^x \W_2(\mu_{t},\mu )^2\Big\} \ge c\sum_{i=1}^\infty \ff 2 {\ll_i^2},\end{equation}
where $c>0$ is a constant which equals to $1$ when $\pp M$ is empty or convex.
Below we prove these two estimates  respectively.

(a)  Let $M$ be compact. Since $\sum_{i=1}^\infty\ff 2 {\ll_i^2}=\infty$ for $d\ge 4$, we only consider   $d\le 3$.
As shown in (b) in the proof of  \eqref{CLT},
we only need to consider  $\nu= h_\nu \mu$ with $\|h_\nu\|_\infty\le C$ for some constant $C>0$. Let $r_t=t^{-\aa}$ for $t\ge 1$ and some $\aa\in (1,2)$. By the triangle inequality of $\W_2$,  for any $\vv>0$ we have
\beq\label{TTU} \W_2(\mu_t,\mu)^2\le (1+\vv) \W_2(\mu_{t,r_t,r_t},\mu)^2+ 2(1+\vv^{-1}) \big\{\W_2(\mu_{t,r_t},\mu_{t,r_t,r_t})^2+ \W_2(\mu_t,\mu_{t,r_t})^2\big\}.\end{equation}
This and Lemma \ref{LN4} yield
 $$ \limsup_{t\to\infty} \sup_{\|h_\nu\|_\infty\le C} \big\{t\E^\nu  \W_2(\mu_t,\mu)^2 \big\}\le (1+\vv)\limsup_{t\to\infty} \sup_{\|h_\nu\|_\infty\le C} \big\{t \E^\nu \W_2(\mu_{t,r_t,r_t},\mu)^2 \big\},\ \ \vv>0. $$   Letting $\vv\downarrow 0$ implies
\beq\label{S11} \limsup_{t\to\infty} \sup_{\|h_\nu\|_\infty\le C} \big\{t\E^\nu  \W_2(\mu_t,\mu)^2 \big\}\le
\limsup_{t\to\infty} \sup_{\|h_\nu\|_\infty\le C} \big\{t \E^\nu \W_2(\mu_{t,r_t,r_t},\mu)^2 \big\},\ \ C>0.\end{equation}
Next, by Lemma \ref{L2.1} and noting that $\ff{\d\mu_{t,r_t,r_t}}{\d \mu} = (1-r_t)f_{t,r_t} +r_t$, for the $p>1$ in Lemma \ref{LN3}, we have
\beq\label{*QR} \beg{split} &\E^\nu [\W_2(\mu_{t,r_t,r_t},\mu)^2]\le (1-r_t)^2\E^\nu \int_M\ff{|\nn L^{-1} (f_{t,r_t}-1)|^2}{\M((1-r_t)f_{t,r_t}+r_t,1)}\d\mu\\
&\le  \E^\nu \int_M \Big\{ |\nn L^{-1} (f_{t,r_t}-1)|^2  +   |\nn L^{-1} (f_{t,r_t}-1)|^2\big| \M((1-r_t)f_{t,r_t}+r_t,1)^{-1}-1\big| \Big\}\d\mu \\
&\le \E^\nu \int_M |\nn L^{-1} (f_{t,r_t}-1)|^2  \d\mu
 + \bigg(\E^\nu \int_M |\nn L^{-1} (f_{t,r_t}-1)|^{2p}\d\mu\bigg)^{\ff 1 p}\\
 &\qquad\qquad\qquad\qquad  \times
 \bigg(\E^\nu\int_M \big| \M((1-r_t)f_{t,r_t}+r_t,1)^{-1}-1\big|^{\ff p{p-1}} \d\mu\bigg)^{\ff {p-1}p}.\end{split}\end{equation}
Combining this with Lemmas \ref{L2.2}, \ref{LN1} and \ref{LN3}, we arrive at
$$ \limsup_{t\to\infty}  \sup_{\|h_\nu\|_\infty\le C} \big\{t \E^\nu [\W_2(\mu_{t,r_t,r_t},\mu)^2]\big\}\le \sum_{i=1}^\infty \ff 2{\ll_i^2},\ \ C>0,$$
which together with \eqref{S11} yields
 \beq\label{LLP}\limsup_{t\to\infty}  \sup_{\|h_\nu\|_\infty\le C} \big\{t \E^\nu [\W_2(\mu_{t},\mu)^2]\big\}\le \sum_{i=1}^\infty \ff 2{\ll_i^2},\ \ C>0.\end{equation}
Then \eqref{E1} holds.

(b)   Let $\pp M$ be either convex or empty. In this case, we have
\beq\label{KPP1}  \W_p(\mu, \nu P_r)^2\le  \e^{2K r} \W_p(\mu,\nu)^2,\ \ r>0, p\in [1,\infty),\end{equation}  where $K\ge 0$ is such that $\Ric_V\ge -K$  holds on $M$,
 and  $(\nu P_r)(f):=\mu(P_r f)$ for $f\in \B_b(M),$   see \cite{RS05} for empty $\pp M$ and \cite[Theorem 3.3.2]{W14} for convex $\pp M$.
Since $\mu_{t,r}=\mu_t P_r$,  \eqref{KPP1} and  \eqref{E2-3} imply
$$ \e^{2K r} \liminf_{t\to\infty}\Big\{t \inf_{x\in M} \E^x\W_2(\mu, \mu_{t})^2\Big\} \ge \liminf_{t\to\infty} \Big\{t \inf_{x\in M} \E^x \W_2(\mu, \mu_{t,r})^2\Big\}
 \ge \sum_{i=1}^\infty \ff 2 {\ll_i^2\e^{2r\ll_i}},\ \ r\in (0,1],$$  which gives \eqref{E2} for $C=1$ by letting $r\to 0.$
 In general, by \cite[Theorem 2.7]{CT17}, there exist constants $C,\ll>0$ such that
\beq\label{TC11} \W_2(\mu, \nu P_r)^2\le  C \e^{\ll r} \W_2(\mu,\nu)^2,\ \ r>0.\end{equation} This together with \eqref{E2-3} yields \eqref{E2}.

 \end{proof}

\beg{proof}[Proof of Theorem \ref{T1.1}(2)]  The boundedness of $\rr$ implies that the weak topology is induced by $\W_2$.
So,
$G:= \{\nu\in\scr P: \W_2(\nu,\mu)^2\in A^\circ\}$ is open  while $F:= \{\nu\in\scr P: \W_2(\nu,\mu)^2\in \bar A\}$ is closed  in $\scr P$ under the  weak topology.
Thus, by Lemma \ref{LDP}, it suffices to prove
\beg{enumerate} \item[(i)] For any set $A\subset [0,\infty)$, $\inf\{I_\mu(\nu): \W_2(\nu,\mu)^2\in A\}= \inf_{r\in A}I(r).$
\item[(ii)] For any $\aa\ge 0$, $\{I\le \aa\}$ is a compact subset of $[0,\infty)$. \end{enumerate} Below we prove these two assertions respectively.

For (i). Let $\tt I(r):= \inf_{\W_2(\nu,\mu)^2=r} I_\mu(\nu).$  We have $\inf\{I_\mu(\nu): \W_2(\nu,\mu)^2\in A\}= \inf_{r\in A}\tt I(r).$ So, it suffices to show that $\tt I(r)$ is increasing in $r\ge 0$, so that
$I(r)=\tt I(r).$ Let $r_1>r\ge 0$ such that $\tt I(r_1)<\infty$. For any $\vv>0$ there exists $\nu=f\mu$ with $\W_2(\nu,\mu)^2=r_1$ such that
$I_\nu(\nu)\le \tt I(r_1)+\vv.$ Consider $\nu_s=s\mu +(1-s)\nu $ for $s\in [0,1].$ Since $\W_2(\nu_s,\mu)$ is continuous in $s$ and $\W_2(\nu_0,\mu)^2= r_1>r, \W_2(\nu_1,\mu)^2=0,$ there exists
$s\in [0,1)$ such that $\W_2(\nu_s,\mu)^2=r$. We have $\ff{\d\nu_s}{\d\mu}= s + (1-s)f$, so that
$$\tt I(r)\le I_\mu(\nu_s)= \mu\Big(\ff{(1-s)^2|\nn f|^2}{4(s+(1-s)f)}\Big) \le \mu(|\nn f^{\ff 1 2}|^2) = I_\mu(\nu)\le \tt I(r_1) +\vv.$$ Letting $\vv\downarrow 0$ we prove $\tt I(r)\le \tt I(r_1)$ as desired.

For (ii). Since $\{I\le\aa\}\subset [0,r_0]$ is bounded, it suffices to prove that $\{I\le \aa\}$ is closed. Let $0\le r_n\to r$ as $n\to\infty$ such that $I(r_n)\le \aa$, it remains to prove $I(r)\le \aa.$
Let $\nu_n=f_n\mu \in \scr P$ such that $\W_2(\nu_n,\mu)^2=r_n$ and
$$I_\mu(\nu_n)=\mu(|\nn f_n^{\ff 1 2}|^2)\le I(r_n)+\ff 1 n\le \aa+\ff 1 n.$$ Then, up to a subsequence, $f_n^{\ff 1 2}\to f^{\ff 1 2}$ in $L^2(\mu)$ for some probability density $f$ with respect to $\mu$.
Thus,  $\nu_n\to \nu:=f\mu$ weakly such that $$\W_2(\nu,\mu)^2 =\lim_{n\to\infty}\W_2(\nu_n,\mu)^2=\lim_{n\to\infty} r_n=r,$$
and $I_\mu(\nu) \le \liminf_{n\to\infty}I_\mu(\nu_n)\le \aa.$ Therefore, $I(r)\le\aa$ as desired.
\end{proof}

\beg{proof}[Proof of Theorem \ref{T1.1}(3)]  As shown in step (b) in the proof of \eqref{CLT}, we only need to prove for $C>0$ that
 \beq\label{CLL} \lim_{t\to \infty}\sup_{\|h_\nu\|_\infty\le C} \big|\P^\nu \big(t  \W_2(\mu_{t},\mu)^2 <a\big)-\nu_0\big((-\infty,a)\big)\big|=0,\ \ a\in \R.\end{equation}
Take $r_t=t^{-\ff 3 2},$  and let
$$\tt\Xi(t):= t\int_M\ff{|\nn L^{-1}(f_{t,r_t}-1)|^2}{\scr M((1-r_t)f_{t,r_t}+r_t,1)} \d\mu,\ \ t>0.$$
 By Lemma \ref{LN1} and Lemma \ref{LN3},    we have
\beg{align*} &\lim_{t\to\infty} \sup_{\|h_\nu\|_\infty\le C} \E^\nu\big|\tt \Xi(t)  -\Xi_{r_t}(t)\big| \le  \lim_{t\to\infty} \sup_{\|h_\nu\|_\infty\le C} \bigg(t^p\E^\nu \int_M |\nn L^{-1}(f_{t,r_t}-1)|^{2p}\d\mu\bigg)^{\ff 1 p}\\
 &\qquad\qquad\qquad \times \bigg(\E^\nu\int_M \big|\scr M((1-r_t)f_{t,r_t}+r_t,1)^{-1}-1\big|^{\ff p{p-1}} \d\mu   \bigg)^{\ff {p-1}p}\\
 &  =0.\end{align*} Combining this with
 \eqref{CLT2}, we obtain
 \beq\label{WM*1} \lim_{t\to\infty} \sup_{\|h_\nu\|_\infty\le C} \big|\P^\nu(\tt\Xi(t)<a) - \nu_0((-\infty,a))\big|=0,\ \ a\in\R.\end{equation}
 By Lemma \ref{L2.1} and \eqref{TTU} we obtain
\beg{align*} &\big\{t\W_2(\mu_t,\mu)^2 - (1+\vv) \tt\Xi(t)\big\}^+\le \big\{t \W_2 (\mu_t,\mu)^2 - (1+\vv) \W_2(\mu_{t,r_t,r_t},\mu)^2\big\}^+ \\
 &\le  2(1+\vv^{-1}) \big\{\W_2(\mu_{t,r_t},\mu_{t,r_t,r_t})^2 + \W_2(\mu_{t,r_t},\mu_t)^2\big\}.\end{align*} Combining this with
   \eqref{LK1}, \eqref{LK2},  \eqref{ACC2} and   $r_t=t^{-\ff 3 2}$, we derive
 \beg{align*}& \lim_{t\to\infty} \sup_{\|h_\nu\|_\infty\le C}  \E^\nu \big\{t\W_2(\mu_t,\mu)^2 - (1+\vv) \tt\Xi(t)\big\}^+  \\
 &\le 2 (1+\vv^{-1}) \lim_{t\to\infty} \sup_{\|h_\nu\|_\infty\le C}  \E^\nu\Big[ t\big\{\W_2(\mu_{t,r_t,r_t},\mu_{t,r_t})^2+ \W_2(\mu_t,\mu_{t,r_t})^2\big\}\Big]=0,\ \ \vv>0.
     \end{align*}Therefore,
  \beq\label{WM*2} \lim_{t\to\infty} \sup_{\|h_\nu\|_\infty\le C}  \P^\nu (t\W_2(\mu_t,\mu)^2 \ge (1+\vv) \tt\Xi(t) +\vv) =0,\ \ \vv>0.\end{equation}
  On the other hand,    \eqref{KPP1}  and \eqref{WM*0}  yield
  \beq\label{W*0} \beg{split} &  \lim_{t\to\infty} \sup_{\|h_\nu\|_\infty\le C}  \P^\nu\big(t\W_2(\mu_t,\mu)^2 \le  \Xi_{r_t}(t) -\vv \big)  \\
  &\le \lim_{t\to\infty} \sup_{\|h_\nu\|_\infty\le C}  \big\{\P^\nu\big(t\e^{-2Kr}\W_2(\mu_{t,r},\mu)^2\le \Xi_r(t)-\vv/2\big) +\P^\nu\big(\Xi_r(t)\le \Xi_{r_t}(t)-\vv/2\big) \big\}\\
  & \le \ff 2 \vv \lim_{t\to\infty} \sup_{\|h_\nu\|_\infty\le C}   \E^\nu\big|\Xi_r(t)- \Xi_{r_t}(t)\big|,\ \ r>0.\end{split}\end{equation}
Since $\sum_{k=1}^\infty \ll_k^{-2}<\infty$, \eqref{SX2} implies
$$\lim_{r,r'\downarrow  0}  \sup_{\|h_\nu\|_\infty\le C,t\ge 1} \E^\nu\big|\Xi_r(t)- \Xi_{r'}(t)\big| \le  \lim_{r,r'\downarrow 0} \sum_{k=1}^\infty \ff{2C}{\ll_k^2} |\e^{-2\ll_k r}-\e^{-\ll_k r'}|  =0. $$  So,  letting $r\to 0$ in \eqref{W*0}  we derive
$$ \lim_{t\to\infty} \sup_{\|h_\nu\|_\infty\le C}  \P^\nu\big(t\W_2(\mu_t,\mu)^2 \le  \Xi_{r_t}(t) -\vv \big)  =0,\ \ \vv>0.$$
Combining this with \eqref{CLT2}, \eqref{WM*1} and \eqref{WM*2}, we prove \eqref{CLL}.
\end{proof}

\section{Proof of Theorem \ref{T1.1'}}  Let $d\ge 4$.
As explained in step (b)  in the proof of \eqref{CLT}, it suffices to prove that for any $C>0$ there exist constant $c>0$ such that
\beq\label{ENN}  \sup_{\|h_\nu\|_\infty\le C}  \E^\nu [\W_2(\mu_t,\mu)^2] \le \beg{cases} c t^{-1}\log (1+t), &\text{if}\ d=4,\\
c t^{-\ff 2 {d-2}}, &\text{if}\ d\ge 5,\end{cases} \end{equation} and
\beq\label{ENN'}  \liminf_{t\to\infty}\Big\{t^{\ff 2{d-2} } \inf_{\|h_\nu\|_\infty\le C}  \{\E^\nu[\W_1(\mu_t,\mu)]\}^2\Big\} >0.\end{equation}

\subsection{Proof  of \eqref{ENN}}
By the  triangle inequality,  Lemma \ref{L2.2} and Lemma \ref{LN4},  we find a constant $c>0$ such that for $r_t \in (0, 1)$,
\begin{align*}
&\E^\nu [\W_2(\mu_{t,r_t,r_t},\mu)^2] \le 2 \E^\nu [\W_2(\mu_{t,r_t,r_t}, \mu_{t,r_t} )^2] + 2 \E^\nu [\W_2(\mu_{t, r_t},\mu)^2] \\
&\le c \bigg( r_t  +   \ff{ \|h_\nu\|_\infty}{t} \sum_{i=1}^\infty \ff1 {\ll_i^2\e^{2r_t\ll_i} } \bigg).
\end{align*}
 Notice that for some constants $c_i>0, i=1,\ldots, 5$ we have
\beg{align*}& \sum_{i=1}^\infty \ff 1{\ll_i^2\e^{2r_t\ll_i}}\le c_1 \sum_{i=1}^\infty i^{-\ff 4 d} \e^{-c_2 r_t i^{\ff 2 d}}
 \le c_3 \int_1^\infty s^{-\ff 4 d} \e^{-c_4 r_t s^{\ff 2 d}}\d s\\
 &=c_3  r_t^{2-\ff d 2} \int_{r_t^{\ff d 2}}^\infty u^{-\ff 4 d} \e^{-c_4 u^{\ff 2 d}}\d u \le c_5\big\{ \log  ( 1 +{r_t}^{-1}  ) 1_{\{d=4\}} +  r_t^{-\ff {d - 4} 2} 1_{\{d\ge 5\}}\big\}.\end{align*}
Taking $r_t = \ff {\log (t + 1)}{t} 1_{\{d=4\}} + t^{- \ff 2 {d-2}} 1_{\{d\ge 5\}}, $ we have
$$\limsup_{t\to \infty} \Big\{ \ff t {\log t}  \sup_{\|h_\nu\|_\infty\le C} \E^\nu [\W_2(\mu_t,\mu)^2]\Big\} < \infty, \ \ d=4,$$
and
$$\lim_{t\to \infty} \Big\{t^{\ff 2 {d-2} } \sup_{\|h_\nu\|_\infty\le C} \E^\nu [\W_2(\mu_t,\mu)^2]\Big\} < \infty,\ \ d\ge 5.$$
Therefore, \eqref{ENN} holds for some constant $c>0$.

\subsection{Proof of \eqref{ENN'}}

In this subsection, instead of $\mu_{t, r}$, we use another method to approximate $\mu_t$. For any $t\ge 1$ and $N\in\mathbb N$, we consider $\mu_{N}:= \ff 1 N\sum_{i=1}^N \dd_{t_i},$ where $t_i:=\ff{(i-1)t}N, 1\le i\le N.$ We write
$$\mu_t= \ff 1 N \sum_{i=1}^N \ff Nt \int_{t_i}^{t_{i+1}}\dd_{X_{s}}\d s.$$
By the convexity of $\W_2^2$, which follows from the Kantorovich dual formula \eqref{KT}, we have
$$\W_2(\mu_N,\mu_t)^2\le \ff 1 N \sum_{i=1}^N  \ff N t \int_{t_i}^{t_{i+1}} \W_2(\dd_{X_{t_i}}, \dd_{X_s})^2\d s = \ff 1 t \sum_{i=1}^N   \int_{t_i}^{t_{i+1}} \rr(X_{t_i},  X_s)^2\d s.$$
Moreover, by \eqref{ACC1} and the Markov property we have
$$\E^\nu [\rr(X_{t_i},  X_s)^2]\le c_1\|P_{t_i} h_\nu\|_\infty (s-t_i)\le c_1C (s-t_i),\ \ s\ge t_i, \|h_\nu\|_\infty\le C.$$
Therefore, there exists a constant $c'>0$ such that
\beq\label{NB1}\sup_{\|h_\nu\|_\infty\le C}  \E^\nu [\W_2(\mu_N,\mu_t)^2]\le \ff{c_1C}t \sum_{i=1}^N \int_{t_i}^{t_{i+1}} (s-t_i)\d s \le \ff{c't} N,\ \  N\in\mathbb N, t\ge 1.\end{equation}

On the other hand, since $M$ is compact possibly with a smooth boundary, and $\mu$ is comparable with the volume measure, we find a constant $c>0$ such that for any $r > 0$,
$$\sup_{x\in M} \mu(B(x,r) )\le c r^d.$$  So, \cite[Proposition 4.2]{RE1} (see also \cite[Corollary 12.14]{RE2}) implies
$$ \{\W_1(\mu_N,\mu)\}^2\ge c_2N^{-\ff 2 d},\ \ N\in\mathbb N, t\ge 1.$$
This together with \eqref{NB1} yields
\beg{align*}   \inf_{\|h_\nu\|_\infty\le C}  \{\E^\nu [\W_1(\mu,\mu_t)]\}^2 &\ge \ff 1 2  \inf_{\|h_\nu\|_\infty\le C}   \{\E^\nu [\W_1(\mu,\mu_N)]\}^2- \sup_{\|h_\nu\|_\infty\le C}   \E^\nu [\W_2(\mu_N,\mu_t)^2]\\
&\ge \ff{c_2}{2N^{\ff 2 d}} -\ff{c't}N,\ \ N\in \mathbb N, t\ge 1.\end{align*}
Taking $N=\sup\{i\in\mathbb N: i\le \aa t^{\ff d{d-2}}\}$ for some $\aa>0$, we derive
$$t^{\ff 2 {d-2}}  \inf_{\|h_\nu\|_\infty\le C}   \{\E^\nu [\W_1(\mu,\mu_t)] \}^2 \ge  \ff {c_2}{2\aa^{\ff 2 d}}- \ff {2c'}\aa $$ for large enough $t\ge 1$. Therefore,
$$\liminf_{t\to\infty} t^{\ff 2 {d-2}}  \inf_{\|h_\nu\|_\infty\le C}  \{ \E^\nu [\W_1(\mu,\mu_t)]\}^2 \ge \sup_{\aa>0} \Big(\ff {c_2}{2\aa^{\ff 2 d}}- \ff {2c'}\aa\Big)>0.$$

\section{Proof of Theorem \ref{TN}}

Recently, a PDE proof of AKT theorem based on Lusin approximation has been developed in \cite{AMB}. In our setting, we adapt a similar strategy.
We first prove the following result.

\beg{prp}\label{PN1} Let $M$ be compact with $d=4$ and $\pp M$ either convex or empty. If for any $C>0$ there exist constants  $\vv,\gg>0$ such that for large enough $t$ we have
\beq\label{PL1}    \{\E^\nu  \W_1(\mu_{t,t^{-\gg}},\mu)\}^2\ge \vv \E^\nu\mu(|\nn (-L)^{-1}(f_{t,t^{-\gg}}-1)|^2),\ \ \nu=h_\nu\mu, \|h_\nu\|_\infty\le C,\end{equation}
then there exists a constant $c>0$ such that
\beq\label{PL2} \liminf_{t\to\infty} t(\log t)^{-1} \inf_{x\in M} \{\E^x \W_1(\mu_t,\mu)\}^2 >0.\end{equation}
\end{prp}
\beg{proof} As   shown in step (b) in the proof of \eqref{CLT}, it suffices to prove that for any constant $C>0$,
\beq\label{PL2'} \liminf_{t\to\infty} t(\log t)^{-1} \inf_{\|h_\nu\|_\infty\le C } \{\E^\nu \W_1(\mu_t,\mu)\}^2 >0.\end{equation}
By taking $r=t^{-\gg}$, we deduce from \eqref{CC} and \eqref{PL1} that
$$\inf_{\|h_\nu\|_\infty\le C} \{\E^\nu  \W_1(\mu_{t,t^{-\gg}},\mu)\}^2\ge \ff{\vv(1-cCt^{-1})} t \sum_{i=1}^\infty
\ff {1}{\ll_i^2\e^{2t^{-\gg} \ll_i}}.$$ Since $\ll_i\sim \ss i$ for $d=4$, combining this with \eqref{KPP1} and \eqref{OPL}, we prove    \eqref{PL2'}.\end{proof}

To verify \eqref{PL1},  we need the following fundamental lemma, where the first assertion is known as  Sard's lemma
(see \cite[p130, Excercise 5.5]{Gry}), and the second
  is called Lusin's approximation which is well-known for $M$ being a bounded open domain in $\R^d$ (see \cite{AF, Liu}). For completeness
we include a simple proof of the second assertion for the present Riemannian setting.

\beg{lem} \label{LMN} Let $M$ be a compact Riemannian manifold, and let $p\in (1,\infty)$.
\beg{enumerate} \item[$(1)$] For any $f\in W^{1,p}(\mu_M)$ and $c\in\R$, we have $\mu_M(\{|\nn f|>0, f=c\})=0.$
\item[$(2)$] There exists a constant $K>0$ such that for any $\aa>0$,
\beq\label{SCD} \inf_{\|\nn u\|_\infty\le\aa} \mu_M( \{ u\ne f \} )\le \ff K {\aa^p} \mu_M(|\nn f|^p),\ \ f\in W^{1,p}(\mu_M).\end{equation}
\end{enumerate}
\end{lem}

\beg{proof} Let $M$ be isoperimetrically embedded into $\R^m$ for some $m>d$, where $d$ is the dimension of $M$. Let $N:=\{N^1,\cdots, N^{m-d}\}$
be   smooth vector fields on $\R^m$ such that for any $\theta\in M$,
\beg{enumerate}\item[(i)] $\<N_\theta^i, N_\theta^j\>=1_{i=j}, \ 1\le i,j\le m-d;$
\item[(ii)] $N_\theta^i$ is orthogonal to the tangent space $TM$ of $M$.\end{enumerate}
Let $B_s:=\{r\in \R^{m-d}: |r|<s\}$ and  $M_{s}:=\{x\in\R^m: {\rm dist}(x,M)<s\}, s>0.$
Since $M$ is compact, there exists $s_0>0$ such that the map
$$M\times B_{s_0}\ni (\theta,r)\mapsto \theta+ \<r, N_\theta\>\in M_{s_0}$$
is diffeomorphic, where $\<r, N_\theta\>:=\sum_{i=1}^{m-d} r_i N_\theta^i, r=(r_1,\cdots, r_{m-d})\in \R^{m-d}.$
 We   simply denote $ \theta+ \<r, N_\theta \>$ by its polar coordinate $(\theta,r).$   Let $\LL$ be the Lebesgue measure on the open set $M_{s_0}\subset \R^m$.
Then there exist a constant $c_1\ge 1$ such that   under the polar coordinate   $(\theta,r)\in M\times B_{s_0}$ we have
\beq\label{VLM0} c_1^{-1} \mu_M(\d\theta)\d r \le  \LL(\d\theta,\d r)\le c_1 \mu_M(\d\theta)\d r,\end{equation}
\beq\label{VLM}     |\nn g|^2  (\theta,r)   \le c_1  \big(|\nn g(\cdot,r)(\theta)|^2+ |\nn g(\theta,\cdot)|^2(r) \big),\end{equation}
\beq\label{DIST} c_1^{-1}  \rr(\theta_1,\theta_2)  \le  |\theta_1-\theta_2| \le c_1 \rr(\theta_1,\theta_2),\ \  \theta_1, \theta_2 \in  M.\end{equation}
Now, for any $f\in W^{1,p}(\mu_M)$, we extend it to $M_{s_0}$ by letting
$$\tt f(x)= f(\theta),\ \ \text{if}\ x=(\theta,r).$$
So, by the Lusin approximation result on the bounded open domain $M_{s_0}$ in $\R^m$, there exists a   constant $c_2>0$ such that for any $\aa>0$, we find  a $(c_1^{-1}\aa)$-Lipschtiz function $\tt u$ on
$M_{s_0}$ such that
$$\LL(\{x\in M_{s_0}: \tt u(x)\ne \tt f(x)\})\le \ff {c_2}{\aa^p} \int_{M_{s_0}}| \nn \tt f|^p\d\LL.$$
Combining this with \eqref{VLM0}, \eqref{VLM}  and  noting that $\tt f(\theta,r)$ does not depend on $r$,  we find a   constant $c_3>0$ such that
$$\int_{B_{s_0}} \mu_M( \{ f\ne \tt u(\cdot,r) \} )\d r\le \ff{c_3}{\aa^p}   \int_{M} |\nn f|^p\d\mu_M.$$
 Thus, there exist  $r\in B_{s_0}$ and a constant $c_4>0$  such that
\beq\label{KLL} \mu_M( \{f\ne \tt u(\cdot, r) \} )\le \ff{c_4}{\aa^p}\mu_M(|\nn f|^p).\end{equation}
Since $|\tt u(x)-\tt u(y)|\le c_1^{-1}\aa|x-y|$, \eqref{DIST} implies
$$|\tt u(\theta_1, r)-\tt u(\theta_2,r)|\le c_1^{-1} \aa |\theta_1-\theta_2|\le  \aa\rr(\theta_1,\theta_2),\ \ \theta_1,\theta_2\in M.$$
Therefore,   \eqref{KLL}  implies \eqref{SCD} for $K=  c_4.$
\end{proof}
Besides the proof of AKT theorem in \cite{AMB} using Riesz transform bounds, a simplified  Fourier analytic technique is developed in
  \cite{BL19}. We will follow the line of   \cite{BL19}  to prove the following result, which together with Proposition \ref{PN1} implies the assertion in  Theorem \ref{TN}.
\def\i{{\rm i}}

\beg{prp} Let   $M=\TO^4$ and $V=0$. Then for any $\gg\in (0,\ff 2 7)$, there exists a constant $c>0$ such that
\beq\label{PL0} \{\E^\nu \W_1(\mu_{t,t^{-\gg}},\mu)\}^2 \ge c \E^\nu \mu(|\nn (-\DD)^{-1}(f_{t,t^{-\gg}}-1)|^2),\ \ t\ge 2 \end{equation} holds for any probability measure $\nu$ on $M$. \end{prp}

\beg{proof} (a) Let $f_t= (-\DD)^{-1}(f_{t,t^{-\gg}}-1).$ By Lemma \ref{LMN}, for any $\aa>0$ there exists an $\aa$-Lipschitz function $u$ such that
$$\mu( \{ u-f_t\ne 0 \} )\le \ff K{\aa^4}\mu(|\nn f|^4)$$ and
$$\int_{\{u-f_t=0\}} |\nn (f_t-u)|^4\d\mu=0.$$
Then there exists a constant $K_1>0$ such that
\beg{align*} &\aa\W_1(\mu_{t,t^{-\gg}},\mu)\ge \mu_{t,t^{-\gg}}(u)-\mu(u) =\int_{\TO^4} u(f_{t,t^{-\gg}}-1)\d\mu= \int_{\TO^4} \<\nn u, \nn f_t\>\d\mu \\
&= \mu(|\nn f_t|^2)- \int_{\TO^4}\<\nn f_t, \nn (f_t-u)\>\d\mu= \mu(|\nn f_t|^2)- \int_{\{f_t\ne u\}}\<\nn f_t, \nn (f_t-u)\>\d\mu\\
&\ge \mu(|\nn f_t|^2) - \mu(1_{\{f_t\ne u\}} |\nn f_t|^2) - \aa \mu(1_{\{f_t\ne u\}}|\nn f_t|)\\
&\ge \mu(|\nn f_t|^2) -  \ss{\mu(|\nn f_t|^4)\mu(\{f_t\ne u\})} - \aa\mu(\{ f_t\ne u \})^{\ff 3 4} \mu(|\nn f_t|^4)^{\ff 1 4}\\
&\ge  \mu(|\nn f_t|^2)- 2K_1\aa^{-2}   \mu(|\nn f_t|^4),\ \ \aa>0.  \end{align*}
Thus,
\beg{align*} \E^\nu [\W_1(\mu_{t,t^{-\gg}},\mu)]\ge   \aa^{-1}\E^\nu [ \mu(|\nn f_t|^2)]-2 K_1\aa^{-3}   \E^\nu[ \mu(|\nn f_t|^4)],\ \ \aa>0.
\end{align*}
If we can find a constant $K_2>0$ such that
\beq\label{NBB} \E^\nu\mu(|\nn f_t|^4)\le K_2\{\E^\nu |\nn f_t|^2)\}^2,\ \ t\ge 2, \end{equation}
then we arrive at
$$\E^\nu \W_1(\mu_{t,t^{-\gg}},\mu)\ge  \aa^{-1} \mu(|\nn f_t|^2) - 2\aa^{-3} K_1 K_2 \mu(|\nn f_t|^2)^{2},\ \ \aa>0.$$
Taking $\aa= N \mu(|\nn f_t|^2)^{-\ff 1 2}$ for large enough $N>1$, we prove \eqref{PL0} for some constant $c>0.$

(b) It remains to prove \eqref{NBB}. To this end,
we identify   $\TO$ with $[0,2\pi)$ by the one-to-one map
$$[0,2\pi)\ni s\mapsto \e^{\i s},$$
where $\i$ is the imaginary unit. In this way, a point in $\TO^4$ is regarded  as a point in $[0,2\pi)^4$, so that
$\{ \e^{\i\<m,\cdot\>}\}_{m\in \Z^4}$ consist of an eigenbasis of $\DD$ in the complex $L^2$-space of $\mu$, where $\mu$ is the normalized volume measure on $\TO^4$. We have
$$f_t:= (-\DD)^{-1}(f_{t,t^{-\gg}}-1) = \sum_{m\in \Z^4\setminus\{0\}} b_m \e^{-\i\<m,\cdot\>},\ \ b_m:= \ff{\e^{-|m|^2 t^{-\gg} }}{|m|^2 t} \int_0^t \e^{\i\<m,X_s\>}\d s.$$
Then
 \beg{align*} &|\nn f_t(x)|^2 = -\sum_{m_1, m_2 \in \Z^4\backslash\{0\}} \langle m_1, m_2 \rangle b_{m_1} b_{m_2} \e^{-\i \langle m_1+m_2, x \rangle},\\
&|\nn f_t(x)|^4   = \sum_{m_1,\cdots, m_4\in \Z^4\backslash\{0\}}  \langle m_1, m_2 \rangle   \langle m_3, m_4\> b_{m_1}b_{m_2}b_{m_3}b_{m_4} \e^{-\i \<m_1+m_2+m_3+m_4, x\>}.\end{align*}
Noting that $\mu(\e^{-\i \<m,\cdot\>})=0$ for $m\ne 0$, we obtain
\beq\label{INT2}  \E^\nu \mu(|\nn f_t|^2) = \sum_{m\in \Z^4\setminus \{0\}} |m|^2 \E^\nu[b_m b_{-m}],\end{equation}
\beq\label{INT4}  \E^\nu \mu(|\nn f_t|^4) = \sum_{(m_1,m_2,m_3,m_4)\in \S}  \<m_1,m_2\>\<m_3, m_4\> \E^\nu[b_{m_1}b_{m_2}b_{m_3}b_{m_4}],\end{equation}
where $\S:=\{(m_1,m_2,m_3,m_4)\in \Z^4\setminus\{0\}: m_1+m_2+m_3+m_4=0\}.$

By the definition of $b_m$,  we obtain
$$\E^\nu[b_m b_{-m}]=
\ff   {\e^{-2|m|^2 t^{-\gg}}} {|m|^2t^2 } \int_{[0,t]^2}  \E^{\nu}  \e^{\i\<m,X_{s_2}-X_{s_1}\>}\d s_1\d s_2.$$
Since  the Markov property and $\E^x\e^{\i\<m,X_s\> } = \e^{-|m|^2s}\e^{\i\<m,x\>}$ imply
\beq\label{CDT}  \E^{\nu}  (\e^{\i\<m,X_{s_2}-X_{s_1}\>}|\F_{s_1\land s_2})= \e^{-|m|^2|s_1-s_2|},\end{equation}
we find a constant $\kk>0$ such that
$$ \E^\nu[b_m b_{-m}]=  \ff   {\e^{-2|m|^2 t^{-\gg}}} {|m|^2t^2 }  \int_{[0,t]^2} \e^{-|m|^2|s_1-s_2|} \d s_1 \d s_2
\ge \ff {\kk \e^{-2|m|^2 t^{-\gg}}}{|m|^6t},\ \ t\ge 2. $$
Combining this with \eqref{INT2} we find a constant $\kk_1>0 $ such that
\beq\label{ESR1}  \beg{split}&\E^\nu \mu(|\nn f_t|^2)  \ge \sum_{m\in \Z^4\setminus \{0\}} \ff {\kk \e^{-2|m|^2 t^{-\gg}}}{|m|^4 t} \ge \ff{\kk_1}t\int_1^\infty \ff{\e^{-2s^2t^{-\gg}}}{s}\d s\\
& \ge \ff {\kk_1} {t\e^2} \int_1^{t^{ \ff \gg 2}} s^{-1}\d s= \ff{\kk_1\gg}{2\e^2}(t^{-1} \log t),\ \ t\ge 2.\end{split}\end{equation}

(c) By \eqref{INT4}, to estimate $\E^\nu\mu(|\nn f_t|^4)$, we  calculate  $\E^\nu[b_{m_1}b_{m_2}b_{m_3}b_{m_4}]$ for $(m_1,m_2,m_3,m_4)\in \S$.
Let $D(t)=\{(s_1, s_2, s_3, s_4) \in [0, t]^4: 0\le s_1 \le s_2 \le s_3 \le s_4 \le t\}$, and let ${\bf S}$ be the set of all the permutations of $\{ 1, 2, 3, 4 \}$.   We have
\beq\label{q1} \beg{split}
&\E^\nu [b_{m_1} b_{m_2}b_{m_3} b_{m_4}]  \\
&= \ff{\e^{-\sum_{p=1}^{4} |m_p|^2 t^{-\gg}}}{t^4\prod_{p=1}^{4} |m_p|^2}   \int_{[0,t]^4}\E^\nu [ \e^{i\langle m_1, X_{s_1}\rangle} \e^{i\langle m_2, X_{s_2}\rangle} \e^{i\langle m_3, X_{s_3}\rangle} \e^{i\langle m_4, X_{s_4}\rangle}] \d s_1 \d s_2 \d s_3 \d s_4  \\
& = \ff{\e^{-\sum_{p=1}^{4} |m_p|^2 t^{-\gg}}}{t^4\prod_{p=1}^{4} |m_p|^2}   \sum_{(i,j,k,l) \in {\bf S}} \int_{D(t)} \E^\nu [ \e^{i\langle m_i, X_{s_1}\rangle} \e^{i\langle m_j, X_{s_2}\rangle} \e^{i\langle m_k, X_{s_3}\rangle} \e^{i\langle m_l, X_{s_4}\rangle}] \d s_1 \d s_2 \d s_3 \d s_4.
\end{split}\end{equation}
Similarly to \eqref{CDT}, we have
\beg{align*}  &\E^\nu \big[\e^{\i\langle m_i, X_{s_1}\rangle} \e^{\i\langle m_j, X_{s_2}\rangle}\e^{\i\langle m_k, X_{s_3}\rangle} \e^{\i\langle m_l, X_{s_4}\rangle}\big]\\
&= \E^\nu\big[\e^{\i\langle m_i, X_{s_1}\rangle} \e^{\i\langle m_j, X_{s_2}\rangle}\e^{\i\langle m_k, X_{s_3}\rangle}
\E^\nu(\e^{\i\langle m_l, X_{s_4}\rangle}|\F_{s_3})\big]\\
&=\e^{-|m_l|^2(s_4-s_3)}\E^\nu\big[\e^{\i\langle m_i, X_{s_1}\rangle} \e^{\i\langle m_j, X_{s_2}\rangle}
 \E^\nu(\e^{\i\<m_k+m_l, X_{s_3}\>}|\F_{s_2})\big]\\
&=\e^{-|m_l|^2(s_4-s_3)-|m_l+m_k|^2(s_3-s_2)} \E^\nu\big[\e^{\i\langle m_i, X_{s_1}\rangle}\E^\nu( \e^{\i\langle m_j+m_k+m_l, X_{s_2}\rangle} |\F_{s_1}) \big]\\
&=  \e^{-|m_l|^2(s_4-s_3)-|m_l+m_k|^2(s_3-s_2)-|m_i|^2(s_2-s_1)},\end{align*} where in the last step we have used $m_i+m_j+m_k+m_l=0.$
Combining this with \eqref{q1} we arrive at
\beq\label{q2}\beg{split}&  \E^\nu [b_{m_1} b_{m_2}b_{m_3} b_{m_4}] \\
&=  \ff{\e^{-\sum_{p=1}^{4} |m_p|^2 t^{-\gg}}}{t^4\prod_{p=1}^{4} |m_p|^2}   \sum_{(i,j,k,l) \in {\bf S}} \int_{D(t)} \e^{-|m_l|^2(s_4-s_3)-|m_l+m_k|^2(s_3-s_2)-|m_i|^2(s_2-s_1)} \d s_1 \d s_2 \d s_3 \d s_4.\end{split}\end{equation}
When $m_l+ m_k = 0$, we have
\beg{align*}
& \int_{D(t)} \e^{-|m_l|^2(s_4 - s_3)} \e^{-|m_l+m_k|^2(s_3 - s_2)} \e^{-|m_i|^2(s_2 - s_1)} \d s_1 \d s_2 \d s_3 \d s_4 \\
& = \int_{0}^{t} \int_{s_1}^{t} \int_{s_2}^{t} \int_{s_3}^{t} \e^{-|m_l|^2(s_4 - s_3)}  \e^{-|m_i|^2(s_2 - s_1)} \d s_4 \d s_3 \d s_2 \d s_1
  \le \ff{t^2}{|m_i|^2|m_l|^2}.
\end{align*}
When  $m_l+ m_k \ne 0$, we have
\beg{align*}
& \int_{D(t)} \e^{-|m_l|^2(s_4 - s_3)} \e^{-|m_l+m_k|^2(s_3 - s_2)} \e^{-|m_i|^2(s_2 - s_1)} \d s_1 \d s_2 \d s_3 \d s_4 \\
& = \int_{0}^{t} \int_{s_1}^{t} \int_{s_2}^{t} \int_{s_3}^{t} \e^{-|m_l|^2(s_4 - s_3)} \e^{-|m_l+m_k|^2(s_3 - s_2)} \e^{-|m_i|^2(s_2 - s_1)} \d s_4 \d s_3 \d s_2 \d s_1\\
& \le \ff{t}{|m_i|^2 |m_l+m_k|^2 |m_l|^2}.
\end{align*}
Therefore,  \eqref{q2} implies
 $$\E^\nu [b_{m_1} b_{m_2}b_{m_3} b_{m_4}] \le \ff{\e^{-\sum_{p=1}^{4} |m_p|^2t^{-\gg}}}{\prod_{p=1}^{4} |m_p|^2} \sum_{(i,j,k,l) \in {\bf S}}\Big\{\ff{t^{-2}1_{\{m_l+m_k=0\}}}{|m_i|^2|m_l|^2}+ \ff{t^{-3}1_{\{m_l+m_k\ne 0\}}}{|m_i|^2 |m_l+m_k|^2 |m_l|^2}\Big\},$$
 so that \eqref{INT4} yields
\beq\label{INT4'}   \E^\nu\mu(|\nn f_t|^4) \le C(I_1+I_2), \ \ t\ge 2 \end{equation} for some constant $C>0$, where
\beg{align*} &I_1:=\ff{1}{t^2} \sum_{a, b \in \Z^4 \backslash \{0\}} \ff{1}{|a|^4 |b|^4} \e^{-2(|a|^2+|b|^2)t^{-\gg}} ,\\
&I_2:= \ff{1}{t^3} \sum_{\substack{m_1,m_2,m_3,m_4 \in \Z^4 \backslash \{0\}\\ m_3+m_4 \ne 0}} \ff{\e^{-\sum_{p=1}^{4} |m_p|^2 t^{-\gg}}}{|m_1|^3 |m_2| |m_3| |m_3+m_4|^2 |m_4|^3} .
  \end{align*}
 Obviously, there exist constants $c_1,c_2>0$ such that
 \beq\label{I_1} I_1      \le \ff{c_1}{t^2} \bigg(\int_1^\infty \ff{-2s^2 t^{-\gg}}{s}\d s\bigg)^2 \le c_2 (t^{-1}\log t)^2,\ \ t\ge 2.\end{equation}
 Similarly,
there exists a constant $c_3>0$ such that
 \beg{align*}   \sum_{m_1 \in \Z^4 \backslash \{0\}} \ff{\e^{-|m_1|^2 t^{-\gg} }}{|m_1|^3} \le     c_3 t^{\ff \gg 2} ,\ \
  \sum_{m_2 \in \Z^4 \backslash \{0\}} \ff{\e^{-|m_2|^2t^{-\gg}}}{|m_2|}  \le c_3t^{\ff{3\gg}2},\ \ t\ge 2.\end{align*}  So,
\beq  \label{I_2} \beg{split} I_2  &  = \ff{1}{t^3} \bigg( \sum_{m_1 \in \Z^4 \backslash \{0\}} \ff{\e^{-|m_1|^2 r}}{|m_1|^3}  \bigg) \bigg(\sum_{m_2 \in \Z^4 \backslash \{0\}} \ff{\e^{-|m_2|^2 r}}{|m_2|} \bigg)
   \sum_{\substack{m_3, m_4 \in \Z^4 \backslash \{0\}\\ m_3+m_4 \ne 0}}  \ff{\e^{-(|m_3|^2+|m_4|^2)r}}{|m_3| |m_3+m_4|^2 |m_4|^3}  \\
 &\le c_3^2 t^{2\gg-3} \sum_{m_4 \in \Z^4 \backslash \{0\}} \ff{\e^{-|m_4|^2 r}}{|m_4|^3}  \sum_{m_3 \in \Z^4 \backslash \{0,\{0,-m_4\}-m_4\}} \ff{\e^{-|m_3|^2 r}}{|m_3| |m_3+m_4|^2}.\end{split}\end{equation}
 Write
\beq\label{S_0} \sum_{m_3 \in \Z^4 \backslash \{0,-m_4\}} \ff{\e^{-|m_3|^2 r}}{|m_3| |m_3+m_4|^2} = S_1 + S_2 +S_3,\end{equation}
where
\beg{align*} & S_1 := \sum_{\substack{m_3 \in \Z^4 \backslash \{0,-m_4\} \\ |m_3|\le \ff{|m_4|}2}}  \ff{\e^{-|m_3|^2 t^{-\gg}}}{|m_3| |m_3+m_4|^2},  \\
 &S_2 := \sum_{\substack{m_3 \in \Z^4 \backslash \{0,-m_4\} \\ \ff{|m_4|}2 <|m_3|\le 2|m_4|}}  \ff{\e^{-|m_3|^2  t^{-\gg}}}{|m_3| |m_3+m_4|^2}, \\
 & S_3 := \sum_{\substack{m_3 \in \Z^4 \backslash \{0,-m_4\}  \\ |m_3| > 2|m_4|}}  \ff{\e^{-|m_3|^2  t^{-\gg}}}{|m_3| |m_3+m_4|^2}.  \end{align*}
Since on the region $\{ m_3 \in \Z^4 \backslash \{0,-m_4\}: |m_3|\le \ff{|m_4|}2 \}$ we have $|m_3 + m_4|^2 \sim  |m_4|^2$, there exists a constant $c_4>0$ such that
\beq\label{S_1}
S_1 \le \ff 4 {|m_4|^2} \sum_{m_3 \in \Z^4 \backslash \{0 \}} \ff{\e^{-|m_3|^2  t^{-\gg}}}{|m_3|} \le \ff{c_4 t^{\ff {3\gg}2} }{|m_4|^2},\ \ t\ge 2.
\end{equation}
Next,  since on the region $\{ m_3 \in \Z^4 \backslash \{0,-m_4\}:  |m_3|> 2|m_4| \}$ it holds  $|m_3 +m_4|^2 \sim  |m_3|^2$ and $|m_3|^2 > \ff{|m_3|^2}2 + 2|m_4|^2$, there exists a constant $c_5>0$ such that
$$
S_3 \le  4 \sum_{\substack{m_3 \in \Z^4 \backslash \{0\}\\ |m_3| > 2|m_4|}} \ff{\e^{-|m_3|^2 t^{-\gg}}}{|m_3|^3} \le 4 \e^{- 2|m_4|^2r } \sum_{m_3 \in \Z^4 \backslash \{0 \}} \ff{\e^{-\ff{|m_3|^2 r}2}}{|m_3|} \le  c_5t^{\ff\gg 2} \e^{- 2 |m_4|^2   t^{-\gg}}.
$$ Noting that  $\e^{-s}\le s^{-1}$ for $s>0$,  this implies
\beq\label{S_3} S_3\le \ff{c_5t^{\ff{3\gg}2}}{|m_4|^2},\ \ t\ge 2.\end{equation}
Finally,   on the region $\{  m_3 \in \Z^4 \backslash \{0,-m_4\}: \ff{|m_4|}2 <|m_3|\le 2|m_4|   \}$ there holds $|m_3| \sim  |m_4|$ and $1 \le  |m_3 + m_4| \le 3|m_4|$, so that there exists a constant $c_6>0$ such that
$$
S_2 \le \ff{2\e^{-\ff{|m_3|^2  t^{-\gg}}4}}{|m_4|} \sum_{1 \le  |m_3 + m_4| \le 3|m_4|} \ff1{|m_3 + m_4|^2}  \le c_6  |m_4|\e^{-\ff{ |m_4|^2t^{-\gg}} 8}.
$$
Using  $\e^{- s} \le c s^{-\ff 3 2}$  for some constant $c>0$ and  all $s > 0$, we find a constant $c_7>0$ such that
$$S_2\le  \ff{c_7t^{\ff{3\gg}2}}{|m_4|^2},\ \ t\ge 2.$$  Combining this with \eqref{S_0}, \eqref{S_1} and \eqref{S_3}, we prove
$$I_2\le c_8t^{\ff {3\gg} 2+2\gg-3}= c_8t^{\ff{7\gg }2-3},\ \ t\ge 2 $$ for some constant $c_8>0$.  Since $\gg\in (0,\ff 2 7)$, substituting this and \eqref{I_1} into \eqref{INT4'} we derive
$$\E^\nu \mu(|\nn f_t|^4)\le c_9 (t^{-1}\log t)^2,\ \ t\ge 2 $$ for some constant $c_9>0$.
  This together with \eqref{ESR1} implies \eqref{NBB} for some constant $K_2>0$, and hence finishes the proof.
\end{proof}
\paragraph{Acknowledgement.} We would like to thank Professor Michel Ledoux and the referees for helpful comments.

\end{document}